\documentclass[10pt]{article}
\usepackage{amstext}
\usepackage{amsmath}
\usepackage{amsfonts}
\usepackage{amssymb}
\usepackage{amsthm}
\usepackage{array}
\usepackage{hyperref}
\usepackage{url}
\usepackage{indentfirst}
\usepackage{tikz}
\usepackage{booktabs}
\usepackage{filecontents}
\usepackage[noend]{algpseudocode}
\usepackage[margin=1.5in]{geometry}    % For reducing margin
\usepackage{algorithm}
\usepackage{graphicx}
\usepackage[utf8]{inputenc}
\usepackage{hyperref}
\usepackage{caption}
\usepackage{algorithm}
\usepackage{algpseudocode}
\usepackage{tabu}

\usetikzlibrary{arrows}
\usetikzlibrary{shapes.geometric}
\usetikzlibrary{positioning,arrows.meta,patterns,shadows}
\usetikzlibrary{graphs,graphs.standard,quotes}
\usetikzlibrary{shapes.arrows,patterns}
\usetikzlibrary{positioning,shapes,calc,decorations.pathmorphing}

\usepackage{enumerate} 
\usepackage[justification=centering]{caption} % Figures caption
\usepackage[export]{adjustbox}
\usepackage{caption, subcaption}

\newcommand{\R}{\mathbb{R}}

\newcommand{\cB}{\mathcal{B}}
\newcommand{\cT}{\mathcal{T}}
\newcommand{\om}{\omega}
\newcommand{\Om}{\Omega}

\algrenewcommand\algorithmicrequire{\textbf{Input:}}

\newtheorem{theorem}{Theorem}[section]
\newtheorem{proposition}{Proposition}[section]
\newtheorem{lemma}{Lemma}[section]

\newtheorem{example}{Example}[section]

\newtheorem{remark}{Remark}[section]
\def\cG{\mathcal G}
\def\cP{\mathcal P}
\usepackage{chngcntr}
\counterwithout{remark}{section}
\counterwithout{figure}{section}
\counterwithout{table}{section}
\counterwithout{theorem}{section}
\counterwithout{lemma}{section}
\counterwithout{proposition}{section}

\usepackage[export]{adjustbox}
\usepackage{caption, subcaption}

\newcommand{\cF}{\mathcal{F}}
\newcommand{\cI}{\mathcal{I}}
\newcommand{\cL}{\mathcal{L}}
\newcommand{\cK}{\mathcal{K}}

\newcommand{\cM}{\mathcal{M}}
\newcommand{\cC}{\mathcal{C}}
\newcommand{\cS}{\mathcal{S}}

\def\cG{\mathcal G}
\def\cP{\mathcal P}
\counterwithout{remark}{section}
\counterwithout{figure}{section}
\counterwithout{table}{section}
\counterwithout{theorem}{section}
\counterwithout{lemma}{section}
\counterwithout{proposition}{section}
\counterwithout{example}{section}
\newcommand\correspondingauthor{\thanks{Corresponding author}}

\title{More on discrete convexity}

\author{Vladimir Gurvich\\
vladimir.gurvich@gmail.com\\
National Research University Higher School of Economics,\\
Moscow, Russian Federation,\\RUTCOR, Rutgers University, Piscataway, NJ, United States;\\
\and
Mariya Naumova\correspondingauthor\\
mnaumova@business.rutgers.edu\\
Rutgers Business School, Rutgers University, Piscataway, NJ, United States\\
}

\begin{document}

\maketitle
\begin{abstract}
In several recent papers some concepts  
of convex analysis were extended to discrete sets. 
This paper is one more step in this direction. 
It is well known that a local minimum 
of a convex function is always its global minimum. 
We study some discrete objects 
that share this property and 
provide several examples of convex families 
related to graphs and to two-person games in normal form.  

{\bf AMS subjects}: 91A05, 94D10, 06E30.

\end{abstract}

\section{Hereditary and Convex Discrete Families}

The similarity between convex functions  
and submodular discrete functions 
is actively studied since 1970s; see for example,   
\cite{BEGK03,Edm70,Fuj05,GLS84,Lov83,Mur03,Orl09,Sha71,Sch00,Sch03}.  
Also in several recent papers some concepts and ideas 
of convex analysis were applied to discrete sets and functions 
\cite{DKL03,FKS07,Kos11,Mur21,Mur21a,MMTT20}. 
In both cases matroids play an important role. 

The present paper is another step in this direction. 
It is well known that each local minimum 
of a convex function is always its global minimum. 
We study some discrete sets that have the same property and  
provide several examples related to graphs and 
two-person games in normal form. 

\medskip 

A {\em partially ordered set 
(poset)}  $(\cP, \succ)$ is defined 
by a set $\cP$ and a binary relation $\succ$. 
The latter is assumed to be 
asymmetric, $P \not\succ P$ for any $P \in \cP$, 
and transitive,
$P \succ P'$ and $P' \succ P''$  imply 
$P \succ P''$  for any $P, P', P'' \in \cP$.
Furthermore, $P'$  is called 
a {\em successor} of  $P$ if $P \succ P'$  and 
an  {\em immediate successor}  of $P$ 
if $P \succ P'$  but  $P  \succ P'' \succ P'$  
holds for no  $P'' \in \cP$. 
Respectively,  $P$  is called an 
{\em (immediate) predecessor} of  $P'$  
if  $P'$  is an (immediate) successor of  $P$.  
Notation $P \succeq P'$  means that either  $P \succ P'$  or  $P = P'$. 

Consider an arbitrary {\em finite} poset $\cP$ 
and its subset (family)  $\cF \subseteq  \cP$.
An element  $F \in \cP$  is called a 
{\em (local) minimum} of $\cF$  if  $F \in \cF$  but  
$F' \not \in \cF$  whenever  $F'$  is an (immediate) successor of  $F$. 
We denote by  $\cM(\cF, \cP, \succ)$  and  by  $\cL\cM(\cF, \cP, \succ)$, 
respectively, the set (class)  of all minima and local minima 
of  $\cF$  in  $(\cP, \succ)$. 
Furthermore, we wave some or all arguments of  $\cM$  and $\cL\cM$  
when they are uniquely determined by the context. 

By the above definitions, containment $\cM \subseteq \cL\cM$ always holds.  

\smallskip 

A family  $\cF \subseteq \cP$  is called:  
\begin{itemize}
\item
{\em convex}  if  $\cM(\cF)  = \cL\cM(\cF)$;  

\item
{\em strongly convex}  if  $\cF$  is convex and 
for any  $F \in \cF$  and  $F' \in \cM(\cF)$  
such that  $F \succ F'$  
there exists an immediate successor  
$P$  of  $F$  
such that  $P \in \cF$  and  $F \succ P \succeq  F'$.  
\item
{\em hereditary} if  $P \in \cF$  
whenever $P$  is a successor of some  $F \in \cF$.  
\item
{\em weakly hereditary}  if  $P \in \cF$ 
whenever $F \in \cF, \; F' \in \cM(\cF)$, 
and  $F \succ P \succeq  F'$.  
 \end{itemize}

\medskip 

In accordance with the above definitions,   
the following implications hold: 

\begin{center}
hereditary $\Rightarrow$  weakly hereditary 
$\Rightarrow$  strongly convex  $\Rightarrow$ convex, 
\end{center}

\noindent 
while all inverse implications fail, 
as we will show in this paper. 
Note that the last two concepts become equivalent 
if  $\cF$  and  $\cP$  have the same unique minimum and 
the first two are equivalent whenever  
$\cM(\cF) \supseteq \cM(\cP)$, or more precisely, 
$\cM(\cF, \cP, \succ) \supseteq \cM(\cP, \cP, \succ).$ 

\begin{remark} 
The last concept can be ``slightly" modified as follows: 
\begin{itemize}
\item 
A family $\cF$ is called {\em very weakly hereditary}  
if  $P \in \cF$  whenever  $F \in \cF$ and 
$F \succ P \succeq  F'$  for some  $F' \in \cM(\cF)$.
\end{itemize}

Then, to the above chain of implications 
we can add the following one: 

\begin{center}
hereditary $\Rightarrow$  weakly hereditary $\Rightarrow$ 
very weakly hereditary  $\Rightarrow$ convex.
\end{center}

Yet, we consider this modification only 
in Subsections
\ref{very-weak-2}, \ref{very-weak-1}, and \ref{very-weak}.
since we have few``natural"  
examples of very weakly but not weakly hereditary families,  
although many ``formal" such examples are not difficult to construct. 
\end{remark} 

In the next three sections we consider several examples 
related to directed and non-directed graphs, 
complete edge-chromatic graphs, and two-person 
game forms and normal form games, respectively. 
We survey known results and obtain several new ones. 

\section{Graphs and digraphs} % , and $d$-graphs} 
\subsection{Definitions and preliminaries}
Given a finite (directed) graph $G$, 
we denote by  $V(G)$  and  $E(G)$  
the sets of its vertices and (directed) edges, respectively. 
Multiple edges are allowed but loops are forbidden. 

A (directed) graph  $G$  is called:  
% \begin{itemize}    \item 
{\em null-graph}  if  $V(G) = \emptyset$ and  
%    \item 
{\em edge-free}  if  $E(G) = \emptyset$. 
% \end{itemize}
The null-graph is unique and edge-free, by definition, 
but not vice versa. 

\smallskip 

We will consider two partial orders:  
related to vertices  $\succ_V$  
and to (directed) edges $\succ_E$. 
In the first case,  $G \succ G'$  if  
$G'$  is an induced subgraph  of  $G$, that is, 
$V(G') \subseteq V(G)$  and  $E(G')$  
consists of all (directed) edges of  $E(G)$ 
whose both ends are in $V(G')$.  
In the second case,  $G \succ G'$  if  
$G'$  is a subgraph of $G$  
defined on the same vertex-set, that is, 
$V(G') = V(G)$  and  $E(G') \subseteq E(G)$.   

\medskip 

Given a graph $G$, which may be directed or not, 
and a set of its (induced) subgraphs  $G_1, \dots, G_n$, 
define a family  $(\cF(G), \succ_E)$ 
(respectively,  $(\cF(G), \succ_V)$)    
that consists of all  subgraphs  $G'$  of  $G$ 
containing as  a (induced) subgraph 
at least one of $G_i, i = 1, \dots, n$. 

\begin{lemma} 
\label{WH}
 Both families are weakly hereditary. 
 Furthermore, $(\cF(G), \succ_V)$ 
(respectively,  $(\cF(G), \succ_E)$)  
is hereditary if and only if  $n=1$ 
and  $G_1$  is the null-graph 
(respectively, the edge-free graph).  
\end{lemma}  
 
\proof 
Consider a subgraph  $G' \in  (\cF(G), \succ_V)$, 
(respectively,  $G' \in(\cF(G), \succ_E)$)  that contains 
(as an induced subgraph)  
$G_i$  for some  $i \in [n] = \{1, \dots, n\}.$
Obviously, the above property is kept  
when we delete a vertex from  $V(G') \setminus V(G_i)$
(respectively, an edge  $e \in E(G') \setminus E(G_i)$) if any. 
Obviously, such a vertex (respectively, an edge) 
exists unless  $G = G_i$. 
Thus, in both cases family  $\cF(G)$  is weakly hereditary. 
Obviously, it is hereditary if and only if  
$G_i$  cannot be reduced. 
\qed 

\medskip 

% The concept of $d$-graph will be introduced later, 
% in Section \ref{d-graphs}.  
% In this case only order  $\succ_V$  can be considered. 

\subsection{Connected graphs}  
A graph  $G$  is called {\em connected} 
if for every two distinct vertices $v, v' \in V(G)$ 
it contains a path connecting  $v$  and  $v'$. 
In particular, the null-graph and the one-vertex graphs are connected, 
since they do not have two distinct vertices. 

Recall that $G'$  is 
a {\em spanning tree} of  $G$  if 
$V(G') = V(G), \; E(G') \subseteq E(G)$, and 
$G'$  is a {\em tree}, that is, connected and has no cycles.

\subsubsection{Order $\succ_V$} 
\label{very-weak-2}

As usual, $\cF = \cF(G)$  denote the family 
of all connected induced subgraphs of a given graph  $G$.
Family $\cF(G)$ is strongly convex for any $G$. 
In other words, every connected graph  
$G$  has a vertex $v \in V = V(G)$  
such that $G[V \setminus \{v\}]$  is connected. 
Indeed, $v$ can be any leaf of a spanning tree of  $G$. 
Furthermore, $\cF(G)$ is hereditary 
if and only if  $G$  is complete.
Otherwise, $\cF(G)$ is not even weakly hereditary. 
 
\begin{example}
%\footnotesize
Consider 2-path  $G = (v_1,v_2),(v_2,v_3)$. 
It is connected, 
but $G[\{V \setminus v_2\}] = G[\{v_1,v_3\}]$  is not.
% by deleting  $v_2$  we obtain 
% a not connected subgraph induced by  $\{v_1, v_3\}$. 
Thus, family $\cF(G)$ is not even very weakly hereditary. 

Let us modify our convention 
and assume now that the null-graph is not connected. 
Then, $\cM = \cM(G) = \cL\cM(G)$ and 
$G' \in \cM$ if and only if  $V(G')$  is a single vertex. 
In this case the 2-path  $(v_1,v_2),(v_2,v_3)$ 
becomes very weakly (but not weakly) hereditary. 
Indeed, the target vertex may be $v_2$ 
but not $v_1$ or $v_3$.  
To obtain a not very weakly hereditary 
family  $\cF(G'')$, consider the 3-path  
$G'' = (v_1,v_2),(v_2,v_3),(v_3,v_4)$. 
In this case, each target vertex can be obtained 
by a vertex-eliminating sequence 
that does not respect connectivity. 
\end{example} 

\subsubsection{Order $\succ_E$}  
Given a connected graph  $G$, 
family  $\cF = \cF(G)$  consists of 
all connected subgraphs $G'$ of $G$  with $V(G') = V(G)$. 

%A graph  $G'$  is called a {\em spanning tree} of  $G$  if 
%$V(G') = V(G), \; E(G') \subseteq E(G)$, and 
%$G'$  is a {\em tree}, that is, connected and has no cycles. 

By definition, all spanning trees of  $G$  are in $\cF$  and, 
obviously, they form class  $\cM = \cL\cM$. 
By Lemma \ref{WH}, family  $\cF$  is weakly hereditary. 

\begin{remark}
Let $G$  be a connected graph   
with weighted edges: $w : E(G) \to \mathbb{R}$.  
It is well known \cite{Bor26a,Bor26b,Kru56}  that 
one can obtain a spanning tree of  $G$ 
of maximal total weight by the greedy algorithm, as follows.
Delete an edge  $e \in E(G)$  such that 
(i) $e$  belongs to a cycle of $G$, or in other words, 
the reduced graph is still connected on  $V(G)$, and 
(ii)  $e$  has a minimal weight among all edges satisfying (i). 
Proceed until such an edge exists. 
\end{remark}

\subsection{Disconnected graphs} 
\label{very-weak-1}
\subsubsection{Order $\succ_V$} 
In this case  $\cF = \cF(G)$  is the family 
of all disconnected induced subgraphs  of a given graph  $G$. 
By convention, the null-graph and one-vertex graphs are connected
Hence, class  $\cM(\cF)$  consists of all subgraphs of $G$  
induced by pairs of non-adjacent vertices. 
In particular, $\cF = \emptyset$  if and only if 
there is no such pair, that is, graph $G$  is complete.     

\begin{proposition} 
For every graph  $G$,  
family  $\cF(G)$  is strongly convex. 
\end{proposition} 

\proof 
Consider a not connected induced subgraph  $G'$ of $G$    
and any pair of non-adjacent vertices  $v', v'' \in V(G')$. 
Then, either $V(G') = \{v',v''\}$, 
in which case  $G' \in \cM(\cF(G))$ is minimal,  or 
we will show that 
there exists a vertex  $v \in V(G') \setminus \{v',v''\}$ 
such that subgraph  $G''$  induced by  $V(G') \setminus \{v\}$ 
is still not connected. 
In other words, $\cF(G)$  is strongly convex. 
Assume that $G'$  is not connected and choose 
$w \in V(G')$ such that $w$ and $v'$  
are from distinct 
connected components of $G'$. 
Then, delete  $v \in V(G') \setminus \{v',v'', w\}$, if any. 
The obtained induced subgraph 
still contains  $v'$  and  $v''$.
Furthermore, it is not connected, 
because it still contains $w$. 
It remains to consider the case when
$V(G') = \{v',v'',w\}$. 
Then, we just delete  $w$  
to obtain a not connected induced subgraph.
\qed 

\begin{proposition} 
For every graph  $G$,  
family  $\cF(G)$  is very weakly hereditary. 
\end{proposition} 

\proof 
Consider a not connected induced subgraph  $G'$ of $G$ 
and choose any two vertices $v', v'' \in V(G')$ 
from distinct connected components of  $G'$. 
Then, obviously, every induced subgraph 
$G''$ of $G$  containing both 
$v'$ and  $v''$  is in $\cF$, that is, not connected. 
% Thus, $\cF(G)$  is very weakly hereditary.
\qed 

\medskip 

However, family $\cF(G)$  is not weakly hereditary for some  $G$. 

\begin{example}
Consider graph  $G$  that consists of 
a 2-path  $(v_1,v_2),(v_2,v_3)$  and an isolated vertex $v_0$. 
This graph is disconnected, that is,  $G \in \cF(G)$,  
but, by deleting  $v_0$,  we obtain a connected graph  
$G' \not\in \cF(G)$. Yet, $v_1, v_3 \in V(G')$  and, hence, 
graph  $G''$  induced by these vertices  
is in  $\cF(G)$, moreover, $G'' \in \cM(\cF(G))$. 
Thus, $\cF(G)$  is not weakly hereditary. 

Note that strong convexity holds for  $\cF(G)$, 
because one can delete  $v_2$  rather than  $v_0$. 
\end{example}

\subsubsection{Order $\succ_E$} 
Given a graph  $G$, 
family  $\cF = \cF(G)$  consists of all disconnected 
graphs  $G'$  such that 
$V(G') = V(G)$  and  $E(G') \subseteq E(G)$. 
Then, obviously, family $\cF$  has a unique minimum: 
$\cM(\cF)$  consists of a unique graph, 
which is the edge-free graph on  $V(G)$. 
Obviously, deleting edges and keeping the vertex-set  
respects the non-connectivity. 
Thus, family $\cF$ is hereditary. 

\subsection{Strongly connected directed graphs}  
A directed graph (digraph)  $G$   
is called {\em strongly connected} (SC) 
if for every two (distinct) vertices of  $v,v' \in V(G)$ 
there is a directed path in  $G$  from  $v$  to  $v'$.

\subsubsection{Order $\succ_V$} 
In this case  $\cF = \cF(G)$  is the family 
of all SC induced subgraphs  of a given digraph  $G$. 
This family is not convex for some digraphs. 

\begin{example}
Consider a digraph  $G$  that consists of 
two directed cycles of length at least 3   
with a unique common vertex. 
It is not difficult to verify 
that  $G \in \cL\cM \setminus \cM$.
%Clearly, $G$  is a locally minimal SC digraph, $G \in \cL\cM(\cF(G))$.   
% Indeed,  $G$  is SC but  we destroy this property 
%by deleting any vertex of  $G$.
%Furthermore, $G \not\in \cM$, since each of two cycles of $G$  is in $\cM$. 
%Thus, $G \in \cL\cM \setminus \cM$. 
% that is, family  $\cF$  is not convex. 
\end{example}

\subsubsection{Order $\succ_E$} 
In this case  $\cF = \cF(G)$  is the family 
of all SC subgraphs  $G'$  of a given digraph  $G$ 
such that  $V(G') = V(G)$  and  $E(G') \subseteq  E(G)$. 
Note that  $\cF(G) = \emptyset$  if and only if  $G$  is not SC.  
Obviously, SC subgraphs form 
a monotone non-decreasing subset of $2^E$. 
In other words, for any subgraphs  $G'$  and $G''$  of  $G$  such that 
$V(G') = V(G'') = V(G)$  and  $E(G') \subseteq E(G'') \subseteq E(G)$, 
we have: $G''$  is SC on  $V(G)$  whenever  $G'$  is. 
Thus, by Lemma \ref{WH}, family  $\cF(G)$  is weakly hereditary 
but not hereditary.

\subsection{Not strongly connected directed graphs}  
\subsubsection{Order $\succ_V$} 
\label{very-weak}
In this case  $\cF = \cF(G)$  is the family 
of all not SC induced subgraphs of a given digraph  $G$. 
It is easily seen that  $\cM(\cF)$  consists of all 
pairs of vertices  $v, v' \in V(G)$  such that 
at least one of two arcs  $(v, v')$  or  $(v', v)$  
is missing in $G$.
There are no such pair in  $G$  
if and only if  $\cF(G) = \emptyset$. 

An induced subgraph $G'$ of $G$ is not SC 
(that is, $G' \in \cF$) if and only if  
there exist two (distinct) vertices  $v,v' \in V(G')$
such that there is no directed path from  $v$ to $v'$  in $G'$.
% Hence, an induced subgraph  $G''$  of  $G'$ 
% is in $\cF$  whenever $v, v' \in V(G'')$.
Furthermore, $G'' \in \cM(\cF)$  if and only if 
$G'' = G[{v,v'}]$ is induced 
by distinct two vertices $v,v' \in V(G)$  such that 
either $(v,v') \not\in E(G)$, or $(v',v) \not\in E(G)$, or both.  

Hence, we can reduce $G'$ to $G[{v,v'}]$ 
deleting its vertices, except $v$  and $v'$, in any order.  

Thus, considered family  $\cF$ is very weakly hereditary  
% NOT !!strongly convex. 
and, hence, convex. 

Yet, obviously, it is not hereditary, since 
an induced subgraph of a not SC digraph can be SC. 
The simplest examples are 
two isolated vertices or one arc. 

Moreover, $\cF$ is not even weakly hereditary. 
Consider, for example, 
a directed 3-cycle and one isolated 
(or pending) vertex $v$. 
This digraph is not SC, but after 
deleting  $v$, we obtain a SC digraph. 
Meanwhile any 2 vertices of the 3-cycle 
induce a not SC digraph. 

\begin{proposition} 
Family $\cF$ is strongly convex.
\end{proposition}

\proof 
Consider a not SC digraph $G'$ 
and any its induced subgraph $G'' \in \cM(\cF)$. 
As we know,  $G'' = G[\{v,v'\}]$ 
for some $v,v' \in V(G')$  such that  
either  $v,v'$, or $v',v$, or both  are not in $E(G')$. 
Since digraph $G'$ is not SG, 
there is a vertex  $w \in V(G')$ such that 
there exists no directed path 
either from $w$ to $v$, or from $v$ to $w$, or both.
%Furthermore, there are $w,w' \in V(G')$  such that 
%$G'$  contains no directed path from  $w$  to $w'$.  
%Note that $v$  and $v'$   as well as $w$ and  $w'$  are distinct, 
%while sets $\{v,v'\}$  and  $\{w,w'\}$  may intersect. 
Let us delete vertices of  $G'$, 
except $v,v'$ and $w$, one by one in any order,  
getting $G''' = G[\{v,v',w\}]$ at the end.
All obtained digraphs remain not CS, 
since they contain  $v$ and  $w$. 
Finally, delete $w$ to obtain  
$G'' = G[\{v,v'\}] \in \cM(\cF)$. 
\qed

\subsubsection{Order $\succ_E$} 
In this case  $\cF = \cF(G)$  is the family 
of all not SC subgraphs  $G'$  of a given digraph  $G$ 
such that $V(G') = V(G)$  and  $E(G') \subseteq E(G)$.  
Obviously, all these subgraphs are not SC 
whenever  $G$  is not SC. 
Thus, family $\cF$  is hereditary. 

\subsection{Ternary graphs} 
A graph is called {\em ternary} 
if it contains no induced cycle of length 
multiple to $3$.

By definition, family  $\cT$  
of ternary graphs is hereditary in order $\succ_V$. 
In contrast, in order  $\succ_E$  this family  
is not even convex, as the following example shows.

\begin{example}
From \cite{CSSS20}  we know that 
``D. Kr{\'a}l asked (unpublished): 
Is it true that in every ternary graph 
(with an edge) there is an edge  $e$ such that 
the graph obtained by deleting $e$ is also ternary? 
This would have implied that all ternary graphs are 3-colourable, 
but has very recently been disproved; 
a counterexample was found by M. Wrochna. 
(Take the disjoint union of a 5-cycle and a
10-cycle, and join each vertex of the 5-cycle 
to two opposite vertices of the 10-cycle, in order.)"
In other words, 
consider the standard model of the Petersen graph, 
with two 5-cycles, as in the Figure 1. 
Then, subdivide every edge of the ``outer"  cycle 
by a vertex and connect it 
with the ``opposite" vertex of the ``inner" 5-cycle.

% Maria Chudnovski, Alex Scott, Paoul Seimour, and  Sophi Spirkl, 
% PROOF OF THE KALAI–MESHULAM CONJECTURE, 
% GIL KALAI and ROY MESHULAM. PICTURES.

\begin{figure*}[h]
\centering
\begin{tikzpicture}[scale = 0.4]

\node[regular polygon, regular polygon sides=10, minimum size=4cm, draw=black, shape border rotate=20, name=x] at (6,6) {};
\node[draw=none, regular polygon, regular polygon sides=5, minimum size=1.5cm, name=y] at (6,6) {};

\foreach \corner in {1,2,...,10}
\node[circle,ball color=black,inner sep=0pt,minimum size=5pt] at (x.corner \corner){};

\foreach \corner in {1,2,...,5}
\node[circle,ball color=black,inner sep=0pt,minimum size=5pt] at (y.corner \corner){};

\draw (y.corner 1) -- (y.corner 3);
\draw (y.corner 1) -- (y.corner 4);
\draw (y.corner 2) -- (y.corner 4);
\draw (y.corner 2) -- (y.corner 5);
\draw (y.corner 3) -- (y.corner 5);

\draw foreach[
  evaluate={\cornerA=int(2*\c);
            \cornerB=int(\cornerA-1);
            \cornerC=ifthenelse(\cornerB+5>10, int(mod(\cornerB+5,10)), int(\cornerB+5));}] \c in {1,...,5}{
              (x.corner \cornerB) -- (y.corner \c) -- (x.corner \cornerC)
              };
\end{tikzpicture}
\caption{Wrochna's example}
\label{f1}
\end{figure*}
\end{example}

% Wrochna's example shows that the family 
% of ternary graphs is not convex in order  $\succ_E$. 

\begin{remark} 
In \cite{CSSS20} it was proven that 
chromatic numbers of all ternary graphs are bounded by a constant. 
Yet, it is much larger than 3. 
\end{remark} 

\begin{remark} 
Consider the skeleton graph of the cube. 
Obviously, an induced 6-cycle appears 
whenever we delete an edge. 
However, this graph itself contains two induced 6-cycles. 

Also, an induced 6-cycle appears whenever we delete 
an edge of icosidodecahedron -  
a polyhedron with twenty triangular faces and twelve pentagonal faces, 
which has 30 identical vertices, with two triangles and 
two pentagons meeting at each, and 60 identical edges, 
each separating a triangle from a pentagon (see Fig. \ref{DaVinci}. 
% $https://en.wikipedia.org/wiki/Icosidodecahedron$).
Yet, this graph itself contains triangles and induced 9-cycles.

\begin{figure}[ht!]
\centering
\includegraphics[scale=0.5]{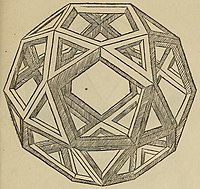}
\caption{Icosidodecahedron. 
Illustration for Luca Pacioli's "Divina proportione" by Leonardo da Vinci \label{DaVinci}}
\end{figure}
\end{remark} 

\subsection{Non-ternary graphs} 

By definition,  non-ternary graph 
contains an induced cycle of length multiple to $3$ 
(a ternary cycle, for short).  
Given a graph  $G$, denote by  
$\cC_3(G)$  (respectively, by $\cI\cC_3(G)$) 
the set of its (induced) ternary cycles 
and by  $\cF(G)$  the family of its  non-ternary subgraphs. 
From Lemma \ref{WH} we will derive that,  
with respect to (wrt) both orders $\succ_V$ and $\succ_E$, 
family  $\cF(G)$  is weakly hereditary but not hereditary. 

\subsubsection{Order $\succ_V$} 
In this case, $\cM(\cF(G)) = \cI\cC_3(G)$. 
Given an induced subgraph  $G'$  of  $G$  
that contains a ternary cycle  $C \in \cI\cC_3(G)$, 
one can delete a vertex  $v \in V(G') \setminus V(C)$  
such that the reduced graph  $G''$  still contains 
$C$  as an induced subgraph unless  $G' = C$. 
This exactly means that family  $\cF(G)$ is weakly hereditary.  
Obviously, it is not hereditary, since deleting a vertex 
might destroy all ternary cycles of  $G'$.  

\subsubsection{Order $\succ_E$} 
In this case, $\cM(\cF(G))$  
is in a one-to-one correspondence with $\cC_3(G)$. 
Recall that   $V(G') = V(G)$  for each subgraph  $G' \in \cF(G)$. 
Hence, $G'$  consists of a cycle $C \in \cC_3(G)$  
and several isolated vertices from  $V(G) \setminus V(C)$. 

Given a non-ternary subgraph  $G''$  of  $G$  
such that  $V(G'') = V(G)$  and  $G''$ contains a 
(not necessarily induced) ternary cycle  $C \in  \cC_3(G)$, 
one can delete an edge $e \in V(G'') \setminus V(C)$  
such that the reduced graph  % $G'''$  
still contains  $C$  unless  $G'' = C$. 
By Lemma \ref{WH}, family  $\cF(G)$ is weakly hereditary. 
Obviously, it is not hereditary, since deleting an edge  
might destroy all cycles in  $G''$ of length multiple of 3.

\subsection{Perfect and imperfect graphs} 
\subsubsection{Definitions and preliminaries} 
Given a graph $G$, as usual,   
$\chi = \chi(G)$ and  $\omega = \omega(G)$  
denote its chromatic and clique numbers, respectively.
Recall that $\chi$ is the minimum number of colors 
in a proper vertex-coloring of $G$ and 
$\omega$ is the number of vertices in a maximum clique of $G$. 
Obviously, $\chi(G) \geq \omega(G)$  for every graph  $G$.  

\medskip 

Graph  $G$  is called {\em perfect}  
if  $\chi(G') = \omega(G')$ 
for every induced subgraph  $G'$  of $G$. %, including  $G$  itself. 

\medskip 

Thus, by definition, in order $\succ_V$ 
the family of perfect graphs is hereditary. 

\medskip 

This concept was introduced in 1961 by Claude Berge  \cite{Ber61} 
(see also \cite{Ber75} for more details)
who made the following two conjectures:

\medskip 

{\bf Perfect Graph Conjecture:} 
$G$  is perfect if (and only if) 
the complementary graph  $\overline{G}$  is perfect. 
It was proven in 1972  by Laslo Lovász  
\cite{Lov72a,Lov72b} and since then is called the 
{\em Perfect Graph Theorem (PGT)}.  

\medskip

A {\em hole} is a cycle of length at least 4 
and an {\em anti-hole} is the complement of a hole.

\medskip 

{\bf Strong perfect graph conjecture:}   
graph  $G$  is perfect if and only if 
it contains no induced odd holes and anti-holes, 
In other words, odd holes and odd anti-holes 
are minimal imperfect graphs in order $\succ_V$. 
This conjecture was proven by 
M. Chudnovsky, N. Robertson, P. Seymour, and R. Thomas in 2002 
and published in 2006 \cite{CRST06}. 
Since then this statement is called the 
{\em Strong Perfect Graph Theorem (SPGT)}.  

\medskip 

A polynomial recognition algorithm for perfect graphs 
was obtained by M. Chudnovsky, G. Cornu{\'e}jols, X. Liu, P. Seymour, and K. Vu{\v s}kovi{\'c} in 2002 and published in 2005 \cite{CCLSV05}.  

\subsubsection{Perfect graphs} 
\paragraph{Order $\succ_V$}  
% By definition, 
This family  $\cF$  is hereditary 
and  $\cM(\cF)$   contains only the null-graph.  

\paragraph{Order $\succ_E$}   
An edge of a perfect graph  $G$ is called {\em critical} 
if deletion of it results in an imperfect graph. 
For example, six edges 
$(v_1, v_2), (v_2, v_3), (v_3, v_4), (v_4, v_5),(v_5, v_1)$, and $(v_1, v_3)$ 
form a perfect graph in which  $(v_1, v_3)$ is a unique critical edge. 
This concept was introduced by Annegret Wagler \cite{Wag01}. 
With Stefan Hougary, she proved that a perfect graph 
has no critical edges if and only if it is {\em Meyniel}, 
that is, every its odd cycle of length 5 or more  
(if any) has at least two chords \cite[Theorem 3.1]{Wag01}. 

\medskip

There are perfect graphs in which all edges are critical. 
Some examples were given in \cite[Figures 2 and 3]{BG09} 
and called {\em Rotterdam} graphs. 
Clearly, these graphs are in  $\cL\cM(\cF)  \setminus \cM(\cF)$  
and, hence, the considered family, 
of perfect graphs in order $\succ_E$, is not convex. 

Furthermore, \cite[ Theorem 4]{BG09} claims 
that every edge of the complement of a Rotterdam graph is critical too. 
In other words, a Rotterdam graph becomes imperfect whenever 
we delete an edge from it or add an edge to it.   

Let us note finally that 
no efficient characterization of 
the non-critical-edge-free perfect graphs is known, 
in contrast to the critical-edge-free ones, which are Meyniel.  
The main result of  \cite{CCLSV05}  
provides a polynomial recognition algorithm for the former family. 

\subsubsection{Imperfect graphs} 
\paragraph{Order $\succ_V$}       
In this case, by the SPGT,  $\cM(\cF)  = \cL\cM(\cF)$    
and this set contains only odd holes and odd anti-holes. 
Again, by Lemma \ref{WH}, $\cF$ is weakly hereditary 
but not hereditary. 

\paragraph{Order $\succ_E$}   
In 1972 Elefterie Olaru 
characterized minimal graphs of this family. 
He proved that it is convex and  $G \in \cM = \cL\cM$    
if and only if  $G$  is an odd hole 
plus $k$ isolated vertices for some $k \geq 0$ 
\cite{Ola72}; see also \cite{BG09,OS84}. 
Thus, by Lemma \ref{WH}, family  $\cF$  is weakly hereditary 
but not hereditary.

Note that the odd anti-holes, except $C_5$, are not in $\cL\cM$, since 
each one has an edge whose elimination would result in 
a graph with an induced odd hole.

\subsubsection{Graphs with $\chi = \omega$} 

It is easily seen that  
$\cM = \cL\cM$  for both orders $\succ_V$  and  $\succ_E$;  
in other words, both families are convex. 

Indeed, for $\succ_V$
(respectively, for $\succ_E$)  
sets $\cM$  and  $\cL\cM$  are equal 
and contain only the null-graph 
(respectively, the edge-free graphs)  
see Propositions 2 and 3 in \cite{BG09}.  

Yet, obviously, deleting a vertex or an edge may fail 
the equality  $\chi = \omega$. 
Thus, both considered families are not hereditary.

\subsubsection{Graphs with $\chi > \omega$} 
\paragraph{Order $\succ_V$} 

By SPGT, every graph with $\chi > \omega$ 
contains an odd hole or odd anti-hole as an induced subgraph; 
in other words,   
class  $\cM$  contains only the odd holes and odd anti-holes. 

Class  $\cL\cM$  is wider; 
it consists of all so-called {\em partitionable} graphs 
defined as follows: 
Graph  $G$  is partitionable if  
$\chi(G)  > \omega(G)$  but 
$\chi(G') = \omega(G')$
for each induced subgraph 
$G'$  of  $G$  such that  
$V(G') = V(G) \setminus \{v\}$  
for a vertex $v \in V(G)$. 
Such definition is one of many 
equivalent characterizations of partitionable graphs;  
this follows easily from the pioneering results of \cite{BHT79,Pad74} 
and it is explicit in \cite{BGH02}. 

Thus, the considered family is not convex. 

\begin{remark} 
The above characterization of $\cM$  is based on SPGT, 
which is very difficult, 
while the case of $\cL\cM$  is simple.
In contrast, partitionable graphs 
are much more sophisticated than the odd holes and anti-holes. 
Although very many equivalent characterizations 
of partitionable graphs are known 
(see, for example, \cite{BHT79,BBGMP98,CGPW84,GT94,Pad74} ) 
yet, their structure is complicated and not well understood. 
For example, the fact that each partitionable graph contains 
an induced odd hole or anti-hole is equivalent with the SPGT. 

The following two questions about partitionable graphs are still open. 
In addition to the odd holes and odd anti-holes, there is one more 
partitionable graph $G_{17}$ that has 17 vertices and has no: 
(i) small transversals and (ii) uncertain edges. 
It is open whether  (i) or (ii) may hold for other partitionable graphs. 
% except odd holes, odd anti-holes, and $G_{17}$
The conjecture that (i) cannot, if true, would significantly strengthen SPGT; 
see \cite{BBGMP98,BGH02,CGPW84,GT94} for the definitions and more details.
\end{remark} 

\paragraph{Order $\succ_E$} 
By Lemma \ref{WH}, the corresponding family  $\cF$  is weakly hereditary:  
class $\cM = \cL\cM$  consists of 
odd holes with $k$ isolated vertices, for some $k \geq 0$. 
This follows from Olaru's Theorem \cite{Ola72};  
see also \cite{OS84} and  \cite[Proposition 1]{BG09}.
Family $\cF$ is not hereditary, since 
obviously, inequality $\chi > \omega$  may 
turn into equality after deleting an edge. 

\subsection{Kernels in digraphs} 
\subsubsection{Definitions and preliminaries}

Given a finite digraph  $G$, 
a vertex-set  $K = K(G) \subseteq V(G)$  is called 
a {\em kernel} of  $G$  if it is (i)  independent and  (ii) dominating, that is, 

(i)  $v, v' \in K(G)$  for no directed edge $(v, v') \in E(G)$   and 

(ii) for every  $v \in V(G) \setminus K(G)$  there is a directed edge 
$(v,v')$  from  $v$  to some  $v' \in K(G)$. 

This definition was introduced in 1901 by Charles Bouton \cite{Bou1901}  
for a special digraph (of the popular game of NIM) and then in 1944 
it was extended by John Von Neumann and Oskar Morgenstern 
for arbitrary digraphs in \cite{NM44}.

It is not difficult to verify 
that an even directed cycle has two kernels, 
while an odd one has none. 
This obvious observation was generalized in 1953 by Richardson \cite{Ric53} as follows:  
A digraph  has a kernel whenever 
all its directed cycles are even. 
The original proof was simplified in \cite{Isb57,HNC65,NL71,GSNL84,BD90,DW12}.
% in 1957 by Isbell \cite{Isb57}, 
% in 1965  by Harary, Norman, and Cartwright  \cite{HNC65}, 
% in 1971  by  Neumann-Lara \cite{NL71}, 
% in 1984  by Galeana Sanchez and Neumann-Lara \cite{GSNL84},  
% in 1990 by Berge and Duchet \cite{BD90}, 
% in 2012 by Dyrkolbotn  and  Walicki \cite{DW12}.  

\begin{remark}
It is not difficult to verify 
that a digraph has at most one kernel 
whenever all its directed cycles are odd \cite{BG06}. 
This claim combined with the Richardson Theorem imply  
that an acyclic digraph has a unique kernel. 
The latter statement is important for game theory, 
allowing to solve finite acyclic graphical zero-sum 
two-person games. 
Of course, it has a much simpler direct proof \cite{NM44}. 
\end{remark}

Already in 1973  V{\'a}sek Chv{\'a}tal proved that 
it is NP-complete to recognize 
whether a digraph has a kernel \cite{Chv73}.  

\subsubsection{Kernell-less digraphs} 
\subsubsection*{Order $\succ_E$} 
In this case, given a digraph  $G$,  
family  $\cF(G)$  contains only the digraphs 
$G'$  such that  $V(G') = V(G), E(G') \subseteq E(G)$, and 
$G'$  has no kernel.
By Richardson's Theorem,  $G' \in \cM(\cF(G))$ 
if and only if  $G'$ is a directed odd cycle in  $G$ 
(plus the set of isolated vertices  $v \in V(G) \setminus  V(G')$).  

In 1980 Pierre Duchet \cite{Duc80} conjectured that 
every kernel-less digraph  $G'$  has an edge  $e \in E(G')$   
such that the reduced digraph  $G'' = G' - e$  
(that is,  $E(G'') = E(G') \setminus \{e\}$) is still kernel-less 
unless  $G''$  is an odd cycle plus  $k$  isolated vertices  
for some $k \geq 0$; 
in other words, family  $\cF$  of kernel-less digraphs is convex, 
$\cM(\cF) = \cL\cM(\cF)$.  
This statement, if true,  
would significantly strengthen Richardson’s theorem. 
Yet, it was shown in \cite{AFG98}  that 
a circulant  with 43 vertices is a counter-example, 
% that is, 
a locally minimal but not minimal kernel-free digraph.

Let us recall that a circulant  
$G = G_n(\ell_1, \dots, \ell_q)$  is defined as a digraph 
with  $n$  vertices, $V(G) = [n] = \{1, \dots , n\}$ 
and $nq$  arcs, 
$E(G) = \{(i, i + j) \mid i \in [n], j \in [q] = \{1, \dots, q\}\}$, 
where standardly all sums are taken $\mod n$.

\begin{example} (\cite{AFG98}) 
% In \cite{AFG98}, 
It was shown that a circulant 
$G_n(1, 7, 8)$  has a kernel if and only if 
$n \equiv 0 \mod 3$ or $n \equiv 0 \mod 29$. 
Hence, $G_{43}(1, 7, 8)$ is kernel-less. 
Yet, a kernel appears whenever an arc of this circulant is deleted. 
Due to circular symmetry, it is sufficient to consider only three cases 
and delete one of the arcs (43, 1), (43, 7), or (43, 8). 
It is not difficult to verify that, respectively, 
the following three subsets become kernels:

\medskip 

$K_1 = \{1, 5, 10, 14, 16, 19, 25, 28, 30, 34, 39, 43\}$,

$K_7 = \{7, 9, 11, 13, 22, 24, 26, 28, 37, 39, 41, 43\}$,

$K_8 = \{3, 5, 8, 14, 17, 19, 23, 28, 32, 34, 37, 43\}   \subseteq \{1, \dots , 43\} = V$. 
\end{example} 

Thus, the set of edge-minimal kernel-free digraphs 
is a proper subset of the locally edge-minimal ones. 
Although, only one digraph from the difference is known, 
it seems that the latter class, 
unlike the former one, is difficult to characterize. 
For example, it is not known whether a circulant $G_n(\ell_1, \ell_2)$ 
can be a locally edge-minimal kernel-less digraph,   
but it is known that it cannot if  $n \leq  1,000,000$ \cite{AFG98}.

\subsubsection*{Order $\succ_V$}  
In this case also family  $\cF$  of the  kernel-less digraphs is not convex. 
Although it seems difficult to characterize 
(or recognize in polynomial time) both classes 
$\cM(\cF)$  and  $\cL\cM(\cF)$  of the 
(locally) vertex-minimal kernel-less digraphs,  
yet, some digraphs from  $\cL\cM \setminus \cM$  
can be easily constructed;  
see for example, \cite{DW12,GS82,GSNL84,GSG07}. 
For completeness we provide one more example. 

\begin{example}
Circulant  $G = G_{16}(1, 7, 8)$ is kernel-less, 
since 16 is not a multiple of 3 or 29. 
Yet, a kernel appears whenever we delete a vertex from  $G$. 
Due to circular symmetry, without loss of generality 
(wlog) we can delete ‘‘the last’’ vertex, 16, and verify that 
vertex-set  $\{1, 3, 5, 7\}$ becomes a kernel.
Hence, $G \in \cL\cM$, but  
$G \not\in \cM$, since $G$ contains a directed triangle, 
$1 + 7 + 8 = 16$, which is kernel-less. 
\end{example}

\subsubsection{Digraphs with kernels} 
We will show that in each order  $\succ_V$  or $\succ_E$  
the corresponding family  $\cF$  
is strongly convex but not weakly hereditary. 

Given an arbitrary digraph  $G$, 
obviously, class  $\cM(\cF)$  contains a unique digraph in both cases:  
the null-graph for $\succ_V$  and 
the edge-free graph with vertex-set $V(G)$ for $\succ_E$. 
Thus, by Lemma \ref{WH}, both families are not weakly hereditary.  

\begin{proposition}
Both families  $(\cF, \succ_V)$  and  $(\cF, \succ_E)$  are strongly convex.    
\end{proposition} 

\proof 
Fix a digraph  $G$  ith a kernel  $K \subseteq V(G)$. 

In order  $\succ_V$  
delete all vertices of  $V(G) \setminus K$, 
if any, one by one. 
By definition, $K$  remains a kernel in every reduced digraph. 
This reduction results in an independent set, 
which is a kernel itself.  
Now we can delete all vertices one by one 
getting the null-graph at the end. 

In order  $\succ_V$, 
first, we delete all arcs within  $V(G) \setminus K$,  
then all arcs from  $V(G) \setminus K$  to  $K$ 
(if any, one by one, in both cases).
By definition, $K$  remains a kernel in every reduced digraph.  
This reduction results in the edge-free digraph 
on the initial vertex-set  $V(G)$. 
\qed 

\subsubsection{On kernel-solvable graphs}
A graph is called {\em kernel-solvable} if 
every its clique-acyclic orientation has a kernel; 
see \cite{BG96,BG98,BG06} for the precise definitions and more details. 
In 1983 Claude Berge and Pierre Duchet conjectured 
that a graph is kernel-solvable if and only if it is perfect.   
The ``only if part" follows easily from the SPGT 
(which remained a conjecture till 2002). 
The "if part"  was proven in \cite{BG96}; see also \cite{BG98,BG06}. 
This proof is independent of SPGT. %, which is very difficult 
The family of perfect graphs is hereditary, by definition.  

As we know, the family of kernel-less digraphs 
is not convex wrt order $\succ_E$,   
in contrast with the family of not kernel-solvable graphs, 
which is convex \cite[Theorem 1]{BG98}.  
% As we know, it is weakly hereditary but not hereditary. 
%  For imperfect graphs this was proven in 1972 by Elefterie Olaru \cite{Ola72}.  

\section{Complete edge-chromatic graphs} 
\label{d-graphs}

\subsection{Definitions and preliminaries}

A $d$-graph $\cG = (V ; E_1, \dots , E_d)$  
is a complete graph whose edges are colored by $d$ colors 
$I = [d] = \{1, \dots, d\}$, or in other
words, are partitioned into $d$  subsets some of which might be empty. 
These subsets are called {\em chromatic components}. 
For example, we call G a 2- or 3-graph if G has only 2, 
respectively, 3,  non-empty chromatic components. 
According to this definition 
order  $\succ_E$  makes no sense for $d$-graphs, 
so we will restrict ourselves by  $\succ_V$.

The following $2$-graph $\Pi$ and $3$-graph $\Delta$    
will play an important role:

$\Pi = (V ; E_1, E_2)$, where 
$V = \{v_1, v_2, v_3, v_4\}, \; 
E_1 = \{(v_1, v_2), (v_2, v_3), (v_3, v_4)\}$, and 
$E_2 = \{(v_2, v_4), (v_4, v_1), (v_1, v_3)\}$;

$\Delta = (V ; E_1, E_2, E_3)$, where 
$V = \{v_1, v_2, v_3\}, \; 
 E_1 = \{(v_1, v_2)\}, E_2 = \{(v_2, v_3)\}$, and $E_3 = \{(v_3, v_1)\}$.

\begin{figure*}[h]
\centering
\begin{tikzpicture}[scale = 0.4]
\tikzset{% This is the style settings for nodes
    cli/.style={circle,ball color=black,inner sep=0pt,minimum size=5pt],draw, general shadow={fill=gray!60,shadow xshift=1pt,shadow yshift=-1pt}},
    c1/.style={very thick,black},
    c2/.style={very thick,red!65!black},
    c3/.style={very thick,green!70!black}}
\node[cli] (w0) at (-6,0) {}; 
\node[cli] (w1) at (-6,-4) {};
\node[cli] (w2) at (-2,-4) {};
\node[cli] (w3) at (-2,0) {};

\node (l0) at (-6.75,0) {$v_2$}; 
\node (l1) at (-6.75,-4) {$v_1$};
\node (l2) at (-1.25,-4) {$v_4$};
\node (l3) at (-1.25,0) {$v_3$};

\draw[c1] (w0) -- (w1);
\draw[c1] (w0) -- (w3);
\draw[c2] (w1) to  node [black,below] {$\Pi$} (w2);
\draw[c1] (w2) -- (w3);
\draw[c2] (w2) -- (w0); 
\draw[c2] (w1) -- (w3); 

\node[cli] (v0) at (5,0) {};
\node[cli] (v1) at (3,-4) {};
\node[cli] (v2) at (7,-4) {};

\node (l00) at (2.25,-4) {$v_1$}; 
\node (l11) at (7.75,-4) {$v_3$};
\node (l22) at (4.25,0) {$v_2$};

\draw[c1] (v0) -- (v1);
\draw[c2] (v1) to  node [black,below] {$\Delta$} (v2);
\draw[c3] (v2) -- (v0); 
\end{tikzpicture}
\caption{$2$-graph $\Pi$ and $3$-graph $\Delta$}
\label{f2}
\end{figure*}

Note that both chromatic components 
of $\Pi$ are isomorphic to $P_4$ and 
that $\Delta$  is a 3-colored triangle.

Let us also remark that, formally, $d$-graphs 
$\Pi(d)$ (respectively, $\Delta(d)$) defined 
for every integer $d \geq 2$  (respectively, $d \geq 3$), 
while $d - 2$ (respectively, $d - 3$) 
of their chromatic components are empty. 
We will omit argument $d$ assuming that it is a fixed parameter.

Both $d$-graphs $\Pi$ and $\Delta$ were first considered, 
for different reasons, in 1967 
by Tibor Gallai in his fundamental paper \cite{Gal67}. 
$\Delta$-free d-graphs are sometimes called {\em Gallai's} graphs. 

The  $\Pi$- and  $\Delta$-free d-graphs have important 
applications to  $d$-person positional games and  
to read-once Boolean functions in case $d=2$   
\cite{BG06,BG09,BGM09,BGM14,GG11,Gur77,Gur78,Gur82,Gur84,Gur09,Gur11}.

\subsection{Substitution for graphs and $d$-graphs}  

Given a graph  $G'$, a vertex $v \in V(G)$, 
and graph  $G''$  such that  
$V(G') \cap V(G'') = \emptyset$, 
substitute  $G''$  for  $v$  in  $G'$  
connecting a vertex $v'' \in V(G'')$  
with a vertex  $v' \in V(G')$  if and only if 
$v$  and  $v'$  were adjacent in  $G'$. 
Denote the obtained graph by  $G = G'(v \to G'')$ 
and call it the {\em substitution} 
of  $G''$  for  $v$  in $G'$  or simply the {\em substitution} 
when arguments are clear from the context. 
% We will also say that  $G = G'(v \to G'')$ 
% is a {\em modular decomposition}  
% of  $G$  via  $G'$  and  $G''$. 

Substitution $\cG = \cG'(v \to \cG'')$  for $d$-graphs  
is defined in a similar way. 
We will see that many important classes of graphs and $d$-graphs 
are closed wrt substitution. 
This will be instrumental in our proofs. 
See Section \ref{modular} for more details.  

\subsection{Complementary connected $d$-graphs} 

We say that a $d$-graph $\cG$  is {\em complementary connected (CC)} if 
the complement of each chromatic component of  $\cG$  is
connected on $V$, in other words, 
if for each two vertices $u, w \in V$  and color 
$i \in [d] = \{1, \dots, d\}$ 
there is a path  between  $u$ and  $w$  without edges of $E_i$.
%  A CC $d$-graph  $\cG$  with  $d \geq 1$  will be called a CC-graph.

By convention, the null-$d$-graph and one-vertex  
$d$-graph are not CC. 
It is easily seen that there is no CC  $d$-graph with two vertices 
and that $\Delta$ (respectively, $\Pi$) is a unique CC $d$-graph 
with three (respectively, four) vertices.

It is also easily seen that 
$\Pi$ and $\Delta$ are minimal CC $d$-graphs, 
that is, they do not contain induced CC subgraphs. 
% In \cite{Gur84}  it was shown that 
Moreover, there are no others. 

\begin{theorem} (\cite{Gur78,Gur84}). 
Every  CC  $d$-graph contains  $\Pi$ or $\Delta$. \qed  
\end{theorem}

\begin{remark} 
In case of  $\Pi$, that is, for  $d=2$,  
the result was obtained earlier 
\cite{CLB81,Sei74,Sum71,Sum73}. 
It was one of the problems on the 1970 Moscow Mathematics Olympiad,  
which was successfully solved by seven high-school students \cite{GT86}.
\end{remark} 

%\bibitem{GT86}
%G. Galperin and A. Tolpygo, Moscow Mathematical Olympiads, 
%A. Kolmogorov, ed.,
% Prosveschenie (Enlightenment), Moscow, USSR, 1986, Problem 72 (in Russian).

%\bibitem{CLB81}
% D. Corneil, H. Lerchsand L. Burlingham, 
% Complement reducible graphs, 
% Discrete Applied Mathematics 3 (1981) 163--174.

This theorem can be strengthened as follows: 

\begin{theorem} (\cite{BG09}).  
Every CC $d$-graph  $\cG$, 
except  $\Pi$  and $\Delta$, contains a vertex  $v$  such that 
the reduced $d$-graph $\cG[V \setminus \{v\}]$  is still CC. \qed 
\end{theorem}

This statement was announced in \cite{Gur78,Gur84}  and 
proven in \cite{BG09}, see also \cite{Gur09}. 
It implies that, by deleting vertices one by one,  
we can reduce every CC $d$-graph 
to a copy of  $\Pi$ or $\Delta$. 

In other words, family  $\cF = \cF^{CC}$ 
of CC $d$-graphs is convex and 
the class  $\cM = \cL\cM$  contains only 
$2$-graph $\Pi$ and $3$-graph $\Delta$. 
Let us show that $\cF_d = \cF^{CC}_d$ is not strongly convex, 
already for  $d=2$. 
Since chromatic components of a $d$-graph may be empty, 
this also shows that  
$F^{CC}_d$  is not strongly convex, for any $d$. 

It is both known and easily seen \cite{Gur09} 
that the family of CC graphs, as well as CC $d$-graphs, 
is closed wrt substitution. 
% see Section \ref{modular} for more details. 
Moreover, $\cG = \cG'(v \to \cG'')$ is CC 
if and only if  $\cG'$  is CC.  

\begin{example}
Consider $2$-graph  $\Pi$  
on the vertex-set  $\{v_1,v_2,v_3,v_4\}$ 
and substitute for  $v_4$  another 2-graph  $\Pi'$   
on the vertex-set  $\{v'_1,v'_2,v'_3,v'_4\}$. 
Obtained $2$-graph  $\cG = \Pi(v_4 \to \Pi')$  is  CC, 
since $\Pi$ and $\Pi'$ are  CC and family  
$\cF_d$  is closed wrt substitution.
Furthermore, it is easy to verify that 
the CC  property disappears 
if we delete  $v_1$, $v_2$  or  $v_3$  from  $\cG$. 
Thus, we cannot reduce  $\cG$  to  $\Pi'$  keeping CC,  
which means that family  $\cF^{CC}_2$  
of the CC $2$-graphs is not strongly convex. 

However, we can reduce  $\cG$  to  $\Pi$  keeping CC, 
in agreement with convexity of  $\cF^{CC}_2$.  
\end{example}

\begin{figure*}[h]
\centering
\begin{tikzpicture}[scale = 0.7]
\tikzset{% This is the style settings for nodes
    cli/.style={circle,ball color=black,inner sep=0pt,minimum size=5pt],draw, general shadow={fill=gray!60,shadow xshift=1pt,shadow yshift=-1pt}},
    c1/.style={very thick,black},
    c2/.style={very thick,red!65!black},
    c3/.style={very thick,green!70!black}}
\node[cli] (v1) at (-2,-2) {};
\node[cli] (v2) at (-2,2) {}; 
\node[cli] (v3) at (2,2) {};

\node (l1) at (-2.75,-2) {$v_1$};
\node (l2) at (-2.75,2) {$v_2$}; 
\node (l3) at (2.75,2) {$v_3$};

\node[cli] (v1p) at (1,-3) {};
\node[cli] (v2p) at (1,-1) {}; 
\node[cli] (v3p) at (3,-1) {};
\node[cli] (v4p) at (3,-3) {};

\node (l1p) at (0.85,-0.25) {$v_2'$};
\node (l2p) at (0.85,-3.75) {$v_1'$}; 
\node (l3p) at (3.15,-0.25) {$v_3'$};
\node (l4p) at (3.15,-3.75) {$v_4'$};

\draw[c1] (v2) -- (v1);
\draw[c1] (v2) -- (v3);
\draw[c2] (v1) -- (v3); 

\draw[c1] (v2p) -- (v1p);
\draw[c1] (v2p) -- (v3p);
\draw[c2] (v1p) -- (v4p);
\draw[c1] (v4p) -- (v3p);
\draw[c2] (v4p) -- (v2p); 
\draw[c2] (v1p) -- (v3p); 

\draw[c2] (v1) -- (v1p);
\draw[c2] (v1) -- (v2p);
\draw[c2] (v1) -- (v3p);
\draw[c2] (v1) -- (v4p);

\draw[c2] (v2) -- (v1p);
\draw[c2] (v2) -- (v2p);
\draw[c2] (v2) -- (v3p);
%\draw[c2] (v2) -- (v4p);

\draw[c2] (v2) to[bend right=90, looseness=2] (v4p);

\draw[c1] (v3) -- (v1p);
\draw[c1] (v3) -- (v2p);
\draw[c1] (v3) -- (v3p);
\draw[c1] (v3) -- (v4p);

\end{tikzpicture}
\caption{$2$-graph $\cG = \Pi(v_4  \to \Pi'$)} 
\label{f3}
\end{figure*}

\subsection{Not CC $d$-graphs} 

Let us denote this family by   $\cF_d = \cF^{not-CC}_d$ 
and show that it is strongly convex but not weakly hereditary if  $d > 1$.  
Of course,  $\cF_1$  is hereditary.    

A two-vertex $d$-graph, that is, 
a single edge, is CC only if  $d = 1$.  

\begin{proposition}
Assume that  $d > 1$  and $|V| \ge 2$. 
Each not CC $d$-graph  $\cG = (V ; E_1, \dots, E_d)$  
contains a vertex $v \in V$ 
such that the sub-$d$-graph
$\cG[V \setminus {v}]$  is still not CC.
\end{proposition}

\proof 
Since  $\cG$ is not CC, there is a color  $i \in [d]$ such that 
graph  $\overline{G_i} = (V , \overline{E_i})$ is not connected. 
As we know, one can eliminate vertices of  $V$  
one by one keeping this property until 
$V$  is reduced to two vertices. 
However, the obtained not CC $d$-graph is still not minimal, 
since by convention, the null-$d$-graph and 
a one-vertex $d$-graph are not CC either.   
Thus, the null-$d$-graph is the only (locally) minimal CC $d$-graph. 
So, by the last two steps, 
we reduce the obtained two-vertex not CC $d$-graph 
to a one-vertex $d$-graph and then 
to the null-$d$-graph. Both are not CC. 
Thus, family  $\cF_d$  is strongly convex for all $d$. 
\qed 

Yet, it is not weakly hereditary whenever  $d>1$. 

\begin{example}
Let us add a vertex  $v_0$  to  $\Pi$  or $\Delta$  
and  connect it to all other vertices by edges of the same color. 
Clearly, the obtain $d$-graph is not CC  if  $d > 1$.
Yet, deleting  vertex $v_0$  from it we obtain 
$\Pi$ or $\Delta$, which are both CC. 
Thus, family $\cF_d$ is not hereditary.
Moreover, it is not weakly hereditary, 
by Lemma \ref{WH} unless  $d =1$.
\end{example}

\subsection{CIS property of $d$-graphs} 
\label{CIS_prop_d}
Given a $d$-graph  $\cG = (V; E_1, \dots, E_d)$, 
choose a maximal independent set $S_i \subseteq V$ 
in every graph $G_i = (V, E_i)$  and denote by  
$\cS = \{S_i \mid i \in [d] = \{1, \dots, d\}\}$ 
the obtained collection; furthermore set  
$S = \bigcap_{i=1}^d S_i$.  
Obviously, $|S|\leq 1$  for every $\cS$.
Indeed, if  $v, v' \in S$  then  
$(v, v') \not \in  E_i$  for all  $i \in [d]$, 
that is, this edge has no color. 

\smallskip 

We say that  $\cG$  has the CIS property and call $\cG$  
a {\em CIS  $d$-graph} if  
$S \neq \emptyset$ for every selection  $\cS$. 
CIS $d$-graphs were introduced in 2006 in \cite{ABG06}; 
see also  \cite{BG09,BGM14,BGZ09,DLZ04,DLZ04a,Gur84,Gur09,Gur11,WZZ09,Zan95}.     

\smallskip 

%\begin{remark}
For  $d=2$, a 2-graph consists of two complementary graphs 
$G_1$ and  $G_2$  on the same vertex-set  $V$. 
In this case  CIS  property means that 
in  $G_i$   every maximal clique  
$C$ intersects every maximal stable set  $S$  for  $i = 1,2$. 
(This explains the name CIS.) 
Obviously, $C$  and  $S$  may have at most one vertex in common.  
% \end{remark}

\subsection{Not CIS  $d$-graphs}

It is easy to verify that $d$-graphs  
$\Pi$  and  $\Delta$  are not CIS, 
while every their sub-$d$-graph is CIS, 
in other words, $\Pi$  and  $\Delta$  are minimal not CIS  $d$-graphs. 
Moreover, they are also locally minimal and there are no other.  

\begin{theorem} (\cite{ABG10}). 
Every not CIS  $d$-graph  
$\cG = (V; E_1 \dots, E_d)$, except  $\Pi$  and  $\Delta$, 
has a vertex  $v \in V$  such that the reduced  $d$-graph 
$\cG[V \setminus \{v\}]$  is  not CIS. 
\qed 
\end{theorem}

In other words, for any $d \ge 2$, 
family $\cF_d = \cF_d^{not-CIS}$  
of not CIS  $d$-graphs is convex and 
class  $\cM(\cF_d) = \cL\cM(\cF_d)$  
consists of only of  $\Pi$ if  $d =2$  and 
of $\Pi$ and $\Delta$ if $d > 2$.   

\medskip 

Interestingly, family  $\cF_d^{CC}$  of CC  
$d$-graphs with $d > 1$ has the same properties: 
it is convex and class  
$\cM(\cF_d^{CC}) = \cL\cM(\cF_d^{CC})$  
contains only  $\Pi$  and  $\Delta$. 
However, these two families differ, 
$F_d = \cF_d^{CC} \neq \cF_d^{not-CIS} = F'_d$. 
Moreover, both differences 
$\cF_d \setminus \cF'_d$  and  $\cF'_d \setminus \cF_d$ 
are not empty, already for $d=2$. 

\begin{example}
Consider the bull-graph 
(also called  A-graph) is self-complementary;  
hence, the corresponding bull $2$-graph  
$\cB$  given on Figure \ref{f4}  is both  CC and CIS. 
Thus,  $\cB \in \cF_2 \setminus \cF'_2$. 

Consider  $2$-graph  $\Pi$  colored by colors 1 and 2;  
add to it a new vertex $v_5$  and connect it 
to four vertices of $\Pi$ 
by four edges of the same color, say 1.   
It is easily seen that 
the obtained $2$-graph $\cG$ is not CC and not CIS.
Thus,  $\cG  \in \cF'_2 \setminus \cF_2$.  

The above two examples also show that 
both set-differences are not empty for every $d \geq 2$, 
since chromatic components may be empty. 

\begin{figure*}[h]
\centering
\begin{tikzpicture}[scale = 0.4]
\tikzset{% This is the style settings for nodes
    cli/.style={circle,ball color=black,inner sep=0pt,minimum size=5pt],draw, general shadow={fill=gray!60,shadow xshift=1pt,shadow yshift=-1pt}},
    c1/.style={very thick,black},
    c2/.style={very thick,red!65!black},
    c3/.style={very thick,green!70!black}}
\node[cli] (w1) at (-16,-1) {}; 
\node[cli] (w2) at (-14,-3) {};
\node[cli] (w3) at (-12,-6) {};
\node[cli] (w4) at (-10,-3) {};
\node[cli] (w5) at (-8,-1) {};
\node[label=below:$\cB$] (lab) at (-12, -7.5) {};

\node (l1) at (-16.75,-1) {$v_1$};
\node (l2) at (-14.75,-3) {$v_2$}; 
\node (l3) at (-12,-7) {$v_5$};
\node (l4) at (-9.25,-3) {$v_3$}; 
\node (l5) at (-7.25,-1) {$v_4$};

\draw[c1] (w1) -- (w2);
\draw[c1] (w2) -- (w3);
\draw[c1] (w3) -- (w4);
\draw[c1] (w4) -- (w5);
\draw[c1] (w2) -- (w4); 

\draw[c2] (w1.west) edge [out=-80, in=-90] (w3.north);
\draw[c2] (w1) -- (w4);
\draw[c2] (w1) -- (w5);
\draw[c2] (w2) -- (w5);
\draw[c2] (w5.east) edge [out=-100, in=-90] (w3.north);

\node[cli] (w0) at (-4,-2.5) {}; 
\node[cli] (w1) at (-4,-6.5) {};
\node[cli] (w2) at (0,-6.5) {};
\node[cli] (w3) at (0,-2.5) {};
\node[cli] (w4) at (-2,-0.5) {};

\node (l1_) at (-4.75,-6.5) {$v_1$};
\node (l2_) at (-4.75,-2.5) {$v_2$}; 
\node (l4_) at (0.75,-2.5) {$v_3$}; 
\node (l5_) at (0.75,-6.5) {$v_4$};
\node (l5_) at (-2,0.25) {$v_5$};

\node[label=below:$\cG$] (lab2) at (-2, -7.5) {};
\draw[c1] (w0) -- (w1);
\draw[c1] (w0) -- (w3);
\draw[c2] (w1) -- (w2);
\draw[c1] (w2) -- (w3);
\draw[c2] (w2) -- (w0); 
\draw[c2] (w1) -- (w3);
\draw[c1] (w4) -- (w0); 
\draw[c1] (w4) -- (w1); 
\draw[c1] (w4) -- (w2); 
\draw[c1] (w4) -- (w3); 
\end{tikzpicture}
\caption{2-graphs $\cB$ and $\cG$}
\label{f4}
\end{figure*}
\end{example}

For each  $d \geq 2$  
both families  $\cF_d$  of CC  $d$-graphs and  
$\cF'_d$  of not CIS  $d$-graphs are convex,  
and $\cF_d$  is not strongly convex, already for  $d=2$. 
It remains only to prove that  $\cF'_2$  not strongly convex. 
% \smallskip 
It was shown in \cite{ABG10}  
that CIS $d$-graphs are closed wrt substitution. 
% see Section \ref{modular} for more details.

\begin{example}
% To show that  $\cF'$  is not strongly convex,  
Consider the bull $2$-graph  $\cB$  defined by the edges 

\smallskip 

$(v_1,v_2),(v_2,v_3),(v_3,v_4),(v_2,v_5),(v_3,v_5)$  of color  $1$,  

\smallskip 

$(v_2,v_4),(v_4,v_1),(v_1,v_3),(v_1,v_5),(v_4,v_5)$  of color  $2$. 

\smallskip 
\noindent
Note that vertices $\{v_1, v_2, v_3, v_4\}$  induce a $\Pi$.  
As we already mentioned, $\cB$ is a CIS  $2$-graph, 
while  $\Pi$ is not.
Let us substitute  $v_5$  in  $\cB$  by $2$-graph  $\Pi'$ 
defined by the edges:    

\smallskip 

$(v'_1,v'_2),(v'_2,v'_3),(v'_3,v'_4)$  of color  $1$,  

\smallskip 

$(v'_2,v'_4),(v'_4,v'_1),(v'_1,v'_3)$  of color  $2$. 

\smallskip 
% $(v_6,v_7),(v_7,v_8),(v_8,v_9))$  of color  $1$,  
% $(v_7,v_9),(v_9,v_6),(v_6,v_8)$  of color  $2$. 

\begin{figure*}[h]
\centering
\begin{tikzpicture}[scale = 0.5]
\tikzset{% This is the style settings for nodes
    cli/.style={circle,ball color=black,inner sep=0pt,minimum size=5pt],draw, general shadow={fill=gray!60,shadow xshift=1pt,shadow yshift=-1pt}},
    c1/.style={very thick,black},
    c2/.style={very thick,red!65!black},
    c3/.style={very thick,green!70!black}}
\node[cli] (w1) at (-13,1) {}; %v4
\node[cli] (w2) at (-6,1) {}; %v3
%\node[cli] (w3) at (-12,-6) {};
\node[cli] (w4) at (-6,6) {}; %v2
\node[cli] (w5) at (-13,6) {}; %v1

\node (l1) at (-13,0.25) {$v_4$};
\node (l2) at (-13,6.75) {$v_1$}; 
\node (l3) at (-6,6.75) {$v_2$};
\node (l4) at (-6,0.25) {$v_3$}; 

\node (l1_) at (0.75,1.7) {$v_4'$};
\node (l2_) at (0.75,5.3) {$v_1'$}; 
\node (l3_) at (3.25,1.7) {$v_3'$};
\node (l4_) at (3.25,5.3) {$v_2'$}; 

\node[cli] (v4p) at (1,2.5) {};
\node[cli] (v1p) at (1,4.5) {}; 
\node[cli] (v2p) at (3,4.5) {};
\node[cli] (v3p) at (3,2.5) {};

\draw[c1] (w1) -- (w2);
\draw[c1] (w2) -- (v1p);
\draw[c1] (w2) -- (v2p);
\draw[c1] (w2) to[bend right=60, looseness=0.5] (v3p);

\draw[c1] (w2) -- (v4p);

\draw[c1] (v1p) -- (w4);
\draw[c1] (w4) to[bend right=-60, looseness=0.5] (v2p);

\draw[c1] (v3p) -- (w4);
\draw[c1] (v4p) -- (w4);

\draw[c1] (w4) -- (w5);
\draw[c1] (w2) -- (w4); 

\draw[c2] (w1.west) -- (v1p.north);
\draw[c2] (w1) to[bend right=67, looseness=1] (v2p);
\draw[c2] (w1) to[bend right=67, looseness=1] (v3p);
\draw[c2] (w1.west) -- (v4p.north);

\draw[c2] (w1) -- (w4);
\draw[c2] (w1) -- (w5);
\draw[c2] (w2) -- (w5);

\draw[c2] (w5.east) -- (v1p.west);
\draw[c2] (w5) to[bend right=-67, looseness=1] (v2p);
\draw[c2] (w5) to[bend right=-67, looseness=1] (v3p);
\draw[c2] (w5.east) -- (v4p.west);

\draw[c1] (v2p) -- (v1p);
\draw[c1] (v2p) -- (v3p);
\draw[c2] (v1p) -- (v4p);
\draw[c1] (v4p) -- (v3p);
\draw[c2] (v4p) -- (v2p); 
\draw[c2] (v1p) -- (v3p); 

\end{tikzpicture}
\caption{$2$-graph $\cB'$}
\label{f6}
\end{figure*}

The resulting $2$-graph $\cB' = \cB(v_5 \to \Pi)$  is not CIS. 
Indeed, it is easily seen that two disjoint vertex-sets  
$C = \{v_2,v_3,v'_2,v'_3\}$  and  $S = \{v_1,v_4,v'_1,v'_4\}$ 
form in $\cB'$  maximal cliques of colors 1 and 2, respectively.   

In contrast, we obtain a CIS  $2$-graph, 
by substituting  $v_5$  in $\cB$  by a  
proper sub-$2$-graph  $\cG$  of $\Pi'$. 
Indeed   $\Pi'$  is minimal not CIS, 
hence, $\cG$  is  CIS,  $\cB$  is CIS too, 
and  CIS $d$-graphs are closed wrt substitution. 

Summarizing, we conclude that $2$-graph  $\cB'$  is not CIS, 
but one obtains a CIS sub-$2$-graph by deleting 
any vertex  $v \in \{v'_1,v'_2,v'_3,v'_4\}$ from $\cB'$. 
Hence, if we want to stay in $\cF'(\cB')$, 
we can only delete a vertex from  $V(\Pi) = \{v_1,v_2,v_3,v_4\}$, 
keeping  $\Pi'$  but destroying  $\Pi$.   
Thus, $\cB'$  cannot be reduced to  $\Pi$  within  $\cF'(\cB')$, 
which means that family  $\cF'_2$  is not strongly convex.  

However, in agreement with convexity of family 
$\cF'_2(\cB')$, one can reduce  $\cB'$  to  $\Pi'$  
staying within this family. % $\cF'_2(\cB')$.  

Similarly, we can substitute  $v_5$ in $\cB$  by $\Delta$, 
on vertices  $v'_1, v'_2, v'_3$  edge-colored arbitrarily. 
In particular, we can use colors 1,  or  2, or any other.   
Again, it is not difficult to verify that 
the resulting  $d$-graph graph  
$\cB'' = \cB(v_5 \to \Delta)$  is not CIS. 
(No CIS $d$-graph with $\Delta$ is known, 
see Section \ref{Delta}.) 
Yet, deleting a vertex of $\Delta$ from $\cB''$ 
% , that is, $v'_1, v'_2$ or $v'_3\}$, 
we obtain a CIS sub-$d$-graph of  $\cB''$. 
Indeed, $\Delta$ is minimal not CIS, 
hence, any its sub-$d$-graph is CIS, 
bull $2$-graph $\cB$ is CIS too, and 
CIS $d$-graphs are closed wrt substitution.
Hence, if we want to stay in  $\cF_3(\cB'')$ then 
$\cB''$  can be reduced to  $\Delta$   
induced by $v'_1, v'_2, v'_3$,  
but not to $\Pi$, 
induced by  $\{v_1, v_2, v_3, v_4\}$, 
This disproves the strong convexity of  $\cF_3(\cB'')$. 

\begin{figure*}[h]
\centering
\begin{tikzpicture}[scale = 0.37]
\tikzset{
    cli/.style={circle,ball color=black,inner sep=0pt,minimum size=5pt],draw, general shadow={fill=gray!60,shadow xshift=1pt,shadow yshift=-1pt}},
    c1/.style={very thick,black},
    c2/.style={very thick,red!65!black},
    c3/.style={very thick,green!70!black},
    c4/.style={very thick,blue!70!black},
    c5/.style={very thick,yellow!70!black}}

\node[cli] (w1) at (-33,5) {};
\node[cli] (w2) at (-29,5) {};
\node[cli] (w3) at (-29,1) {};
\node[cli] (w4) at (-33,1) {}; 

\node (l1) at (-33,5.75) {$v_1$}; 
\node (l2) at (-29,5.75) {$v_2$};
\node (l3) at (-29,0.25) {$v_3$}; 
\node (l4) at (-33,0.25) {$v_4$};

\node[cli] (v1p) at (-26,3.75) {}; 
\node[cli] (v2p) at (-24,3) {};
\node[cli] (v3p) at (-26,2.25) {};

\node (l1_) at (-26,4.5) {$v_1'$}; 
\node (l2_) at (-23.25,3) {$v_2'$}; 
\node (l3_) at (-26,1.5) {$v_3'$};

\draw[c1] (w1) -- (w2);
\draw[c1] (w2) -- (w3); 
\draw[c1] (w3) -- (w4); 
\draw[c2] (w1) -- (w4); 
\draw[c2] (w1) -- (w3); 
\draw[c2] (w2) -- (w4);

\draw[c3] (v1p) -- (v3p);
\draw[c1] (v1p) -- (v2p);
\draw[c2] (v2p) -- (v3p);

\draw[c1] (w2) -- (v1p);
\draw[c1] (w2) to[bend right=-55, looseness=1] (v2p);
\draw[c1] (w2) -- (v3p);

\draw[c1] (w3) -- (v1p);
\draw[c1] (w3) to[bend right=55, looseness=1] (v2p);
\draw[c1] (w3) -- (v3p);

\draw[c2] (w1) -- (v1p);
\draw[c2] (w1) to[bend right=-67, looseness=1] (v2p);
\draw[c2] (w1) -- (v3p);

\draw[c2] (w4) -- (v1p);
\draw[c2] (w4) to[bend right=67, looseness=1] (v2p);
\draw[c2] (w4) -- (v3p);
%%%%%%---------------------------
\node[cli] (w1_) at (-20,5) {};
\node[cli] (w2_) at (-16,5) {};
\node[cli] (w3_) at (-16,1) {};
\node[cli] (w4_) at (-20,1) {}; 

\node (l1_) at (-20,5.75) {$v_1$}; 
\node (l2_) at (-16,5.75) {$v_2$};
\node (l3_) at (-16,0.25) {$v_3$}; 
\node (l4_) at (-20,0.25) {$v_4$};

\node[cli] (v1p_) at (-13,3.75) {}; 
\node[cli] (v2p_) at (-11,3) {};
\node[cli] (v3p_) at (-13,2.25) {};

\node (l1_1) at (-13,4.5) {$v_1'$}; 
\node (l2_2) at (-10.25,3) {$v_2'$}; 
\node (l3_3) at (-13,1.5) {$v_3'$};

\draw[c1] (w1_) -- (w2_);
\draw[c1] (w2_) -- (w3_); 
\draw[c1] (w3_) -- (w4_); 
\draw[c2] (w1_) -- (w4_); 
\draw[c2] (w1_) -- (w3_); 
\draw[c2] (w2_) -- (w4_);

\draw[c1] (v1p_) -- (v3p_);
\draw[c3] (v1p_) -- (v2p_);
\draw[c4] (v2p_) -- (v3p_);

\draw[c1] (w2_) -- (v1p_);
\draw[c1] (w2_) to[bend right=-55, looseness=1] (v2p_);
\draw[c1] (w2_) -- (v3p_);

\draw[c1] (w3_) -- (v1p_);
\draw[c1] (w3_) to[bend right=55, looseness=1] (v2p_);
\draw[c1] (w3_) -- (v3p_);

\draw[c2] (w1_) -- (v1p_);
\draw[c2] (w1_) to[bend right=-67, looseness=1] (v2p_);
\draw[c2] (w1_) -- (v3p_);

\draw[c2] (w4_) -- (v1p_);
\draw[c2] (w4_) to[bend right=67, looseness=1] (v2p_);
\draw[c2] (w4_) -- (v3p_);

%%%%%%---------------------------
\node[cli] (w1_2) at (-7,5) {};
\node[cli] (w2_2) at (-3,5) {};
\node[cli] (w3_2) at (-3,1) {};
\node[cli] (w4_2) at (-7,1) {}; 

\node (l1_2) at (-7,5.75) {$v_1$}; 
\node (l2_2) at (-3,5.75) {$v_2$};
\node (l3_2) at (-3,0.25) {$v_3$}; 
\node (l4_2) at (-7,0.25) {$v_4$};

\node[cli] (v1p_2) at (0,3.75) {}; 
\node[cli] (v2p_2) at (2,3) {};
\node[cli] (v3p_2) at (0,2.25) {};

\node (l1_1_) at (0,4.5) {$v_1'$}; 
\node (l2_2_) at (2.75,3) {$v_2'$}; 
\node (l3_3_) at (0,1.5) {$v_3'$};

\draw[c1] (w1_2) -- (w2_2);
\draw[c1] (w2_2) -- (w3_2); 
\draw[c1] (w3_2) -- (w4_2); 
\draw[c2] (w1_2) -- (w4_2); 
\draw[c2] (w1_2) -- (w3_2); 
\draw[c2] (w2_2) -- (w4_2);

\draw[c4] (v1p_2) -- (v3p_2);
\draw[c3] (v1p_2) -- (v2p_2);
\draw[c5] (v2p_2) -- (v3p_2);

\draw[c1] (w2_2) -- (v1p_2);
\draw[c1] (w2_2) to[bend right=-55, looseness=1] (v2p_2);
\draw[c1] (w2_2) -- (v3p_2);

\draw[c1] (w3_2) -- (v1p_2);
\draw[c1] (w3_2) to[bend right=55, looseness=1] (v2p_2);
\draw[c1] (w3_2) -- (v3p_2);

\draw[c2] (w1_2) -- (v1p_2);
\draw[c2] (w1_2) to[bend right=-67, looseness=1] (v2p_2);
\draw[c2] (w1_2) -- (v3p_2);

\draw[c2] (w4_2) -- (v1p_2);
\draw[c2] (w4_2) to[bend right=67, looseness=1] (v2p_2);
\draw[c2] (w4_2) -- (v3p_2);

\end{tikzpicture}
\caption{Three versions of $3$-graph $\cB''$}
\label{f7}
\end{figure*}
\end{example}

\subsection{$\Pi$- and $\Delta$-free $d$-graphs}

Clearly, the family of 
$\Pi$- and $\Delta$-free $d$-graphs is hereditary; 
class  $\cM = \cL\cM$  contains only null-$d$-graph. 

\medskip 

We know that  for two different convex families,  
CC  and not CIS  $d$-graphs, 
class  $\cM = \cL\cM$  contains 
only $d$-graphs  $\Pi$  and $\Delta$. 
Hence, the following three properties 
of a $d$-graph  $\cG$ are equivalent: 

\medskip 

(i)   $\cG$  is  $\Pi$- and $\Delta$-free; 

\smallskip 

(ii)  $\cG$  contains no  CC sub-$d$-graph; 

\smallskip 

(iii) $\cG$  contains only CIS  subgraphs. 

\medskip 

This result allows us to construct one-to-one correspondences between  

\medskip 

(j)   $\Pi$- and $\Delta$-free $d$-graphs; 

\smallskip 

(jj)  vertex  $d$-colored rooted trees; 

\smallskip 

(jjj) tight and rectangular game forms of $d$ players.

\medskip 

For examples, see \cite[Figures 1-6]{Gur09},   
and  \cite[Figures 11-13]{ABG06}. 
In its turn, this result enables us to characterize   
read-once Boolean functions, when  $d=2$,  
\cite{GG11,Gur77,Gur78,Gur84,Gur91}  
and normal forms of graphical 
$d$-person games modelled by trees 
\cite{ABG06,Gur78,Gur82,Gur84,Gur09,Gur11}.  

\subsection{More on CIS $d$-graphs}

Several examples of  CIS $d$-graphs can be found in 
\cite[Section 1.1, Figures 1,2,6]{ABG06}. 
For example, each  $\Pi$- and $\Delta$-free  $d$-graph is CIS.
(Furthermore, according to $\Delta$-conjecture, 
no CIS $d$-graph contains $\Delta$; see Section \ref{Delta} below.) 

We have no efficient characterization 
or recognition algorithm for CIS  $d$-graphs, even for  $d=2$.
The problem looks difficult because 
family  $\cF_2^{CIS}$  of CIS $2$-graphs is not hereditary. 
For example, bull $2$-graph is CIS, but deleting 
its ``top vertex"  
we obtain sub-$2$-graph  $\Pi$, which is not CIS.
Moreover, family $\cF_2^{CIS}$  is not even convex. 

\begin{example}
For any integer $n > 2$, we will construct a $2$-graph $\cG_n$  
such that  $\cG_n \in \cL\cM(\cF_2^{CIS}) \setminus \cM(\cF_2^{CIS})$. 
To do so, consider 
complete bipartite  $n \times n$  graph  $K_{n,n}$, 
its line graph  $G_n = L(K_{n,n})$, 
and its complement $\overline{G_n}$. 
These two complementary graphs on the same vertex-set 
form the required $2$-graph  
$\cG_n = \cL(\cK_{(n,n)}$  for each  $n > 2$; 
see \cite[Figure 1.4]{ABG06} as an example for $n = 3$. 

It is easy to verify that $\cG_n$ is a CIS $2$-graph,
but for each $v \in V(\cG_n)$  
the reduced sub-2-graph $\cG_n[V \setminus {v}]$  is not CIS. 
Due to symmetry, it is enough to check 
this claim for just one arbitrary $v \in V(\cG_n)$.  
Thus, $L(\cG_n)$  is a locally minimal CIS $2$-graph. 
Yet, it is not minimal, since only the null-$2$-graph is.
(Moreover, every $2$-graph with at most $3$ vertices is CIS and  
the only minimal not CIS  $2$-graph is $\Pi$.)  
\end{example}

No efficient characterization of 
locally minimal CIS  $d$-graphs is known. 
However, $\Delta$-conjecture, if true, 
would allow us to reduce arbitrary $d$  to  $d=2$; 
see subsection \ref{Delta} below. 

\medskip 

Let us note finally that a $2$-graph  $\cG = (V; E_1, E_2)$   
is CIS whenever each maximal clique 
of its chromatic component  $G_i = (V,E_i)$ 
has a simplicial vertex, for $i=1$  or  $i=2$ \cite{ABG06}. 
Hence, every  $2$-graph  $\cG$  
is a subgraph of a CIS $2$-graph  $\cG'$. 
This is easy to verify   
\cite[Proposition 1 and Corollary 1]{ABG06}. 
Note, however, that the size of $\cG'$  
is exponential in the size of  $\cG$. 

Thus, CIS $d$-graphs cannot be described 
in terms of forbidden subgraphs, already for $d=2$. 
This is not surprising, since this family  % of CIS $d$-graphs 
is not hereditary. 

\medskip
 
Given a CIS $d$-graph  $\cG$  and a partition $\cP$
of its colors $[d] = \{1, \dots, d\}$  
into  $\delta$  non-empty subsets such that $2 \leq \delta \leq d$,   
merging colors in each of this subsets we obtain a $\delta$-graph  
$\cG' = \cG'(\cG,\cP)$, which we will call 
the {\em projection} of $\cG$  wrt color-merging $\cP$.    
It is not difficult to verify 
(see \cite{ABG06}  for details) that:  

\begin{itemize}
    \item 
    if  $\cG$  is CIS then  $\cG'$  is, but not vice versa; 
    \item 
    if $\cG'$  contains a $\Delta$  then  $G$  does, 
but not vice versa; 

    \item 
    if  $\cG'$  contains a $\Pi$  then  
$\cG$  contains a $\Pi$ or $\Delta$.  
\end{itemize}

We can reformulate the first two claims as follows: 
$\cG'$ is CIS or, respectively, Gallai's  
whenever  $\cG$  is.

\medskip 

For $\delta = 2$ the first claim implies that  
merging an arbitrary set of chromatic components of 
a CIS $d$-graph results in a CIS graph. 

\subsection{Modular decomposition of Gallai's  
($\Delta$-free) $d$-graphs} 
\label{modular}

The operation of substitution $G = G'(v \to G'')$ for graphs  
and $\cG = \cG'(v \to \cG'')$  $d$-graphs was already defined above. 
It can be similarly introduced for multi-variable functions 
and for many other objects; see \cite{Mor85} for more details; 
In this paper $G''$  and  $\cG''$  
are referred to as {\em modules} 
and substitution as {\em modular decomposition}.

\medskip 

We say that a family  $\cF$  of graphs or digraphs 
is exactly closed wrt substitution if  
$G \in \cF$   if and only if  $G', G'' \in \cF$ and, 
respectively, 
$\cG \in \cF$  if and only if  $\cG', \cG'' \in \cF$.
It is both known and easy to verify 
that the following families 
are exactly closed wrt substitution: 
perfect, CIS, and  $P_4$-free graphs;   
CIS, CC, Gallai ($\Delta$-free), $\Pi$- and $\Delta$-free  $d$-graphs. 

Recall that $\cG = \cG'(v \to  \cG'')$ is CC 
if and only if $\cG'$ is CC. 

\medskip 

Recall that a graph  $G$  is called {\em CIS} 
if  $C \cap S \neq \emptyset$  for every 
maximal clique  $C$  and 
maximal stable set  $S$  of  $G$. 
Given a CIS (respectively, CC) 
$2$-graph  $\cG = (V; E_1, E_2)$, 
each of its two chromatic components 
$G_i = (V,E_i), \; i = 1,2$  
is a CIS (respectively, CC)  graph. 
See more about CIS and CC graphs in 
\cite{ABG06,ABG10,BG09,BGM14,BGZ09,DLZ04,DLZ04a,Gur84,Gur09,Gur11,WZZ09,Zan95}. 
   
\medskip

Perfect graphs are closed wrt complementation, 
by the Perfect Graph Theorem \cite{Lov72a,Lov72b}. 
Obviously, CIS, CC, and $P_4$-free graphs 
also have this property, just by definition.  

\medskip 

Let  $\cF$  be a family  of Gallai $d$-graphs 
$\cG = (V; E_1, \dots, E_d)$  such that 
the family $\cF'$  of their chromatic components 
$G_i = (V, E_i), \; i \in [d] = \{1, \dots, d\}$  is 
(i)  closed wrt complementation and 
(ii)  exactly closed wrt substitution. 
For example, family $\cF'$ that contains only  
perfect, or CIS, or CC, or $P_4$-free graphs 
has properties (i) and (ii). 

\medskip 

Every $d$-graphs  $\cG \in \cF$  with  $d \geq 3$  
can be decomposed by two non-trivial $d$-graphs, that is, 
$\cG = \cG'(v \to \cG'')$, where  $d$-graphs  $\cG'$ and $\cG''$ 
are distinct from  $\cG$  and 
from the trivial one-vertex $d$-graph. 

As a corollary, we conclude that Gallai's $d$-graphs 
whose chromatic components have properties (i) and (ii) 
can be decomposed by the $2$-colored such $d$-graphs.  
In other words, for Gallai's $d$-graphs, 
the case of arbitrary  $d$  can be reduced to  $d=2$. 
This statement follows from the results of  
Gyárfás and Simonyi \cite{GS04}, 
which, in their turn, are based on the results of  
Cameron, Edmonds, and Lovasz \cite{CE97,CEL86}, 
Möring  \cite{Mor85}, and Gallai \cite{Gal67}; 
see more details in \cite{ABG06,ABG10} and \cite[Section 4]{Gur09}. 

\medskip 

For example, the modular decomposition 
of $\Pi$- and $\Delta$-free $d$-graphs
has important applications in theory of 
positional (graphical)  $n$-person  games 
modelled by trees;  
in particular, it is instrumental in characterizing 
the normal forms of these games  
\cite{ABG06}, \cite[Remark 3]{Gur75}, 
\cite[Chapter 5]{Gur78}, 
and \cite{Gur82,Gur84,Gur09}.   

\subsection{$\Delta$-conjecture} 
\label{Delta}

All CIS $d$-graphs are Gallai's, or in other words,
$d$-graphs containing  $\Delta$  are not CIS.

\medskip 

This conjecture was suggested in 
\cite[remark on page 71, after the proof of Claim 17]{Gur78}   

and it remains open. 
Some partial results were obtained in  \cite{Gur78}. 
In particular: 

\medskip 

(i) Every $\Pi$- and $\Delta$-free $d$-graph is CIS.

\medskip 

(ii) It is sufficient 
to prove $\Delta$  conjecture for 3-graphs;  
then, it follows for $d$-graphs with arbitrary  $d$. 

\medskip

The second claim was proven by Andrey Gol'berg 
(private communications) in 1975 as follows:  
Consider the projection 
$\cG' = \cG'(\cG,\cP)$ of $\cG$  wrt a color-merging $\cP$.     
As we know, $\cG'$  is CIS  whenever  $G$  is. 
Suppose $\Delta$-conjecture fails for $\cG$, 
in other words, $\cG$  is CIS but not Gallai, 
that is, it contains a $\Delta$, say $\Delta_0$. 
Consider a color-merging $\cP$ with $\delta = 3$  such that 
three colors of  $\Delta_0$  are still pairwise distinct in  $\cP$. 
Then, projection  $\cG' = \cG'(\cG, \cP)$ still contains a $Delta$, 
but  $\cG'$  is a CIS 3-graph. 
Thus, $\Delta$-conjecture fails for 3-graphs too.

\medskip

According to the previous subsection,  
each CIS  $d$-graph is a modular decomposition 
(that is, a superposition of % successive 
substitutions) 
of CIS 2-graphs, modulo $\Delta$-conjecture.
If it holds, studying CIS $d$-graphs 
is reduced to studying CIS graphs.  
Yet, the latter is still difficult. 

\medskip 

Let us note that in case of perfect, CC or $P_4$-free 
chromatic components of a $d$-graph, 
we still have to require that it is $\Delta$-free; 
yet, in case of CIS this requirement may 
be waved, modulo $\Delta$-conjecture; 
see more details in \cite{ABG06,BGZ09,DLZ04,DLZ04a,Gur09,Gur11}.   

\section{Finite two-person normal form games and game forms} 
In this section we consider matrices, 
that is, mappings  $M : I \times J \rightarrow R$, 
whose rows   
$I = \{i_1, \dots, i_n\}$  and columns 
$J = \{j_1, \dots, j_m\}$  are the strategies 
of Alice and Bob, respectively, while 
$R$  may vary: 
it is real numbers $\R$ or their pairs $\R^2$  
in case of matrix and bimatrix games, respectively, 
and  $R = \Om = \{\om_1, \dots, \om_k\}$   
is a finite set of outcomes in case of game forms. 
In all cases, $\cP = \cP(M) = I \cup J$  
is the ground set and 
$succ$  is the containment order over $\cP$; 
in other words, 
$\cP(M)$  consists of all submatrices of  $M$. 
By convention, we identify all elements of  $\cP$  
with $I = \emptyset$  or  $J = \emptyset$:  
they correspond to the empty submatrix, 
which is the unique minimum in  $(\cP, \succ)$. 

\subsection{Saddle points and Nash equilibria}

\subsubsection{Saddle point free matrices} 

In this case $R = \R$  is the set of real numbers.  
We assume that Alice is the maximizer and 
she controls the rows, while 
Bob is the minimizer and he controls the columns. 

An entry of  $M$ is a {\em saddle point} 
(SP) if and only it is 
minimal in its row and maximal in its column   
(not necessarily strictly, in both cases). 
It is well known that a matrix has a SP 
if and only if its maxmin and minmax are equal. 
Obviously, a $2 \times 2$  matrix $M$  has no 
SP if and only if one of its diagonals is 
strictly larger than the other, that is, 
$[r_{i_1,j_1}; r_{i_2,j_2}] \cap 
 [r_{i_1,j_2}; r_{i_2,j_1}] = \emptyset$. 

In 1964 Lloyd Shapley  \cite{Sha64}  proved that 
a matrix has a SP  if 
(but not only if)  every its $2 \times 2$ 
submatrix has a SP. 
In other words, for the family  
$\cF$  of all SP free matrices,  
class  $\cM(\cF)$  of the minimal SP free matrices 
consists of the $2 \times 2$  SP free matrices.    
This result was strengthen as follows: 

\begin{theorem} (\cite{BGM09}).  
Every SP free matrix of size larger that $2 \times 2$ 
has a row or column such that 
it can be deleted and the remaining matrix is still SP free.
\qed 
\end{theorem}

In other words, $\cL\cM(\cF) = \cM(\cF)$, that is, 
family  $\cF$  is convex. 
Yet, it is not strongly convex, as the following example shows.

\begin{example}
Consider the following $4 \times 4$  0,1-matrix $M$:  

\medskip 
\begin{equation*}
\begin{bmatrix}
0 & 1 & 0 & 0\\
1 & 0 & 0 & 0\\
1 & 1 & 0 & 1\\
1 & 1 & 1 & 0\\
\end{bmatrix}
\end{equation*}

\medskip 

The following observations are easy to verify:   

Matrix  $M$ is  SP free and it contains two 
(locally) minimal SP free $2 \times 2$ sabmatrices: 
the first one,  $M_1$, upper left,  
is determined by the first two rows and columns of  $M$, 
while the second one, $M_2$, lower right, 
is determined by the last two rows and columns of  $M$.  

Furthermore, 
if we eliminate one of the last 
(respectively, first) 
two rows or columns of $M$,
then a SP appears 
(respectively, it does not) 
in the obtained submatrix. 
Hence, % eliminating successively rows and columns and 
keeping SP freeness  % of the obtained submatrices, 
one can reduce  $M$  to $M_2$  but not to  $M_1$. 
The first claim is in agreement with convexity, 
while the second one disproves strong convexity 
of the family of SP free matrices.

Moreover, no SP appears 
when we delete the first two rows and columns 
from  $M$ in an arbitrary order, thus 
reducing  $M$  to  $M_2$. 
In other words, family $\cF(M)$ is very weakly hereditary, 
not only convex. 
We leave open, if this hold for property $\cF$ in general. 
\end{example}

\subsubsection{Matrices with saddle points} 

\begin{proposition} 
Given a matrix  $M$  with a SP,  
family $\cF = \cF(M)$  
of all submatrices of  $M$  with a SP 
is strongly convex and very weakly hereditary. 
\end{proposition} 

\proof  Obviously, class  $\cM = \cM(\cF(M))$  
consists of all  $1 \times 1$  submatrices 
(that is, entries) of  $M$. 
By assumption $M \in \cF$. 
Let  $(i^*, j^*)$  be a SP in  $M$.  
Obviously, it remains a SP when we 
delete from  $M$  a row distinct from $i^*$ 
or a column distinct from $j^*$. 
Thus, family $\cF$ is very weakly hereditary and, 
hence, it is convex, $\cM = \cL\cM$.  

To prove strong convexity, 
fix an arbitrary entry  $(i_0, j_0)$  in  $M$  and 
reduce $M$  deleting successively its rows, 
except $i_0, i^*$, and columns, except $j_0, j^*$.  
As we already mentioned, $(i^*,j^*)$  remains a SP. 
Consider three cases: 

\smallskip 

If   $i_0 = i^*$  and  $j_0 = j^*$,   
we arrive to the $1 \times 1$  submatrix  
$(i_0, j_0) = (i^*,j^*)$. 

\smallskip 

If  $i_0 \neq i^*$  and  $j_0 \neq j^*$ 
we arrive to the  $2 \times 2$ submatrix 
formed by these two rows and columns.
Note that  $(i^*,j^*)$  is still a SP. 
Then delete row  $i^*$  getting 
a  $1 \times 2$  submarix. 
Clearly, it has a SP, 
which may be not  $(i^*,j^*)$, yet. 
Finally, we delete column $j^*$  getting  $(i_0, j_0)$. 
 
\smallskip 

If   $i_0 = i^*$  or  $j_0 = j^*$  but not both,  
then we arrive to a  
$1 \times 2$ or $2 \times 1$  submatrix 
that consists of  $(i_0, j_0)$  and  $(i^*, j^*)$. 
Note that  $(i^*,j^*)$  is still a SP. 
Then we delete it getting $(i_0, j_0)$ in one step.
\qed 

\medskip

However, family  $\cF(M)$ is not weakly hereditary. 

\begin{example} Consider matrix 
\begin{equation*}
\begin{bmatrix}
0 & 1 & 0\\
0 & 0 & 1
\end{bmatrix}
\end{equation*}
\medskip 
\noindent 
It has two saddle points, both in the first column, 
but, by deleting this column,   
we obtain a SP free matrix. 
 \end{example} 

\subsubsection{Absolutely determined matrices} 
A matrix is called {\em absolutely determined} 
if every its submatrix has a SP.
By Shapley's theorem \cite{Sha64}, 
it happens if and only if 
each $2 \times 2$  submatrix has a SP. 
This condition can by simplified 
in case of symmetric matrices \cite{GL89,GL90,GL92}.
By definition the considered family is hereditary. 
% https://www.academia.edu/15589693/
% Absolutely_determined_matrices

\subsubsection{Nash equilibria free bimatrices} 
A bimatrix game $(A, B)$ is defined as a pair of mappings 
$a : I \times J \rightarrow \R$  and  
$b : I \times J \rightarrow \R$ 
that specify the utility (or payoff) functions 
of Alice and Bob, respectively. 
Now both players are maximizers.

A situation $(i, j) \in I \times J$ 
is called a {\em Nash equilibrium} (NE) 
if no player can improve the result by choosing another strategy
provided the opponent keeps the same strategy, 
that is, if
$a(i, j) \geq a(i', j) \forall i' \in I$  and 
$b(i, j) \geq b(i, j')  \forall j' \in  J$. 

In other words, $i$  is a best response to  $j$  for Alice and 
$j$  is a best response to  $j$ for Bob.

Clearly, Nash equilibria generalize saddle points, 
which correspond to the zero-sum case: 
$a(i, j) + b(i, j) = 0$ for all  $i \in I$ and $j \in J$. 

However, unlike SP free games, 
the minimal NE-free bimatrix games may be larger than
$2 \times 2$. Let us recall an example from \cite{GL90}.
Consider a $3 \times 3$ bimatrix game $(A, B)$ such that

\begin{center}
$b(i_1,j_1) > b(i_1, j_2) \geq b(i_1, j_3)$,\\
$b(i_2, j_3) > b(i_2, j_1) \geq b(i_2, j_2)$,\\
$b(i_3, j_2) > b(i_3, j_3) \geq b(i_3, j_1)$;\\
\vspace{3mm}
$a(i_2, j_1) > a(i_1, j_1) \geq a(i_3, j_1)$,\\
$a(i_1, j_2) > a(i_3, j_2) \geq a(i_2, j_2)$,\\
$a(i_3, j_3) > a(i_2, j_3) \geq a(i_1, j_3).$ 
\end{center}

Naturally, for situations in the same row 
(respectively, column) the values of $b$ 
(respectively, $a$) are compared, since 
Alice  controls rows and has utility function $a$, 
while Bob controls columns and has utility function $b$.

It is easy to verify that:

\smallskip 

$b(i_1, j_1)$ is the unique maximum in row $i_1$ and 

$a(i_1, j_1)$ is a second largest in the column ${j_1}$.

% Similarly,

\smallskip 

$b(i_2, j_3)$ is the unique maximum in row $i_2$ and 

$a(i_2, j_3)$ is a second largest in column $j_3$;

\medskip 

$b(i_3, j_2)$ is the unique maximum in row $i_3$ and 

$a(i_3, j_2)$ is a second largest in column $j_2$;

\smallskip 

$a(i_2, j_1)$ is the unique maximum in column $j_1$ and 

$b(i_2, j_1)$  is a second largest in row $i_2$;

\smallskip 

$a(i_1, j_2)$ is the unique maximum in column $j_2$ and 

$b(i_1, j_2)$ is a second largest in row $i_1$;

\smallskip 

$a(i_3, j_3)$ is the unique maximum in column $j_3$ and 

$b(i_3, j_3)$ is a second largest in row $i_3$.

\medskip  

Consequently, this game is NE-free, 
since no situation is simultaneously the best 
in its row wrt $b$ and in its
column wrt  $a$. 
Yet, if we delete a row or column then a NE appears. 
For example, let us delete $i_1$.
Then the situation $(i_3, j_2)$ becomes a NE. 
Indeed, $b(i_3, j_2)$ is the largest in the row $i_3$ and 
$a(i_3, j_2)$ is a second largest in the column $j_2$, 
yet, the largest, $a(i_1, j_2)$, was deleted. 
Similarly, situations 
$(i_1, j_1), (i_2, j_3), (i_1, j_2), (i_3, j_3), (i_2, j_1)$ 
become NE after deleting lines $i_2, i_3, j_1, j_2, j_3$, respectively. 

Thus, $(A, B)$ is a locally minimal NE-free bimatrix game. 
Moreover, it is also minimal. 
Indeed, one can easily verify that all $2 \times 2$ 
subgames of $(A, B)$ have a NE and, of course, 
$1 \times 2, \; 2 \times 1$, and $1 \times 1$  games always have it.

In general, the following criterion of local minimality holds: 

\begin{theorem} (\cite[Theorem 3]{BGM09}) 
A bimatrix game $(A, B)$ is a locally minimal NE-free game 
if and only if 
it satisfies the following four conditions:

\medskip 

(i) it is square, that is, $|I| = |J| = k$;

\smallskip 

(ii) there exist two one-to-one mappings 
(permutations) 
$\sigma : I \rightarrow J$ and  $\delta  : J \rightarrow I$  
such that their graphs, $gr(\sigma)$ and $gr(\delta)$, are
disjoint in $I \times J$, or in other words, if 
$(i, \sigma(i)) \neq (\delta(j), j)$ for all $i \in I$ and $j \in J$;

\smallskip

(iiia) % the entry  
$a(\delta(j), j)$ is the unique maximum 
in column $j$ and a second largest 
(though not necessarily unique) in row $\delta(j)$;

\smallskip

(iiib) % the entry 
$b(i, \sigma(i))$ is the unique maximum 
in row $i$ and a second largest 
(though not necessarily unique) in column $\sigma(i)$.
\qed 
\end{theorem}

Thus, we have a simple explicit characterization 
of the class  $\cL\cM$  
of the locally minimal NE-free bimatrix games. 
However, not each such game is minimal. 
Indeed, % the strict improvement cycle is formed 
mappings $\sigma$  and $delta$  define 
$2k$  entries of a  $k \times k$  bimatrix 
locally minimal NE-free game. 
Yet, it may contain some smaller 
$k' \times k'$  NE-free subgames; 
see \cite[Section 3]{BGM09} for more details. 

Thus, the family of NE-free bibatrix games is not convex. 
Interestingly, in this case class  $\cM$  
is more complicated than  $\cL\cM$. 
In contrast to the latter, 
no good characterization of the former is known.

\medskip

It is not difficult to verify that 
in the zero-sum case $k$ cannot be larger than 2.  

\medskip 

In contrast, bimatrix games with NE 
are similar to matrix games with SP. 
This family is weakly hereditary but not hereditary. 
The proof from the previous subsection 
can be applied in this case. 

\subsection{Tightness} 
% \subsubsection{Preliminaries and definitions} 
Given a finite set of outcomes $\Om$,  
let  $X$   and $Y$  be finite sets of strategies 
of Alice and Bob, respectively. 
A mapping  $g : X \times Y \rightarrow \Om$  
is called a {\em game form}. 
One can view a game form as a game without payoffs, 
which are not given yet. 

A game form is called {\em tight}  
if its rows and columns form dual hypergraphs. 
Several equivalent definitions of tightness 
can be found, for example, in \cite{GN21,GN22}. 
Nine examples of game forms are given 
in Figure \ref{f1}; 
the first six are tight, the last three are not.

\begin{figure*}[t]
%\centering
\captionsetup[subtable]{position = below}
\captionsetup[table]{position=top}
\hspace*{4em}
\begin{subtable}{0.3\linewidth}
\centering
\begin{tabular}{|c|c|}
\hline
$\om_1$ & $\om_1$ \\ \hline
$\om_2$ & $\om_3$ \\ \hline
\end{tabular}
\caption*{$g_1$}               \label{tab:g1}
\end{subtable}
\hspace*{-3em}
\begin{subtable}{0.3\linewidth}
\centering
\begin{tabular}{|c|c|c|c|}
\hline
$\om_1$ & $\om_1$ & $\om_2$ & $\om_2$ \\ \hline
$\om_3$ & $\om_4$ & $\om_3$ & $\om_4$ \\
\hline
\end{tabular}
\caption*{$g_2$}
\label{tab:g2}
\end{subtable}
\hspace*{-2.5em}
\begin{subtable}{0.3\linewidth}
\centering
\begin{tabular}{|c|c|c|}
\hline
$\om_1$ & $\om_1$ & $\om_3$\\ \hline
$\om_1$ & $\om_2$ & $\om_2$\\ \hline
$\om_3$ & $\om_2$ & $\om_3$ \\ \hline
\end{tabular}
\caption*{$g_3$}               \label{tab:g3}
\end{subtable}

\vspace{0.5cm}

\hspace*{4em}
\begin{subtable}{0.3\linewidth}
\centering
\begin{tabular}{|c|c|c|}
\hline
$\om_1$ & $\om_1$ & $\om_3$\\ \hline
$\om_1$ & $\om_1$ & $\om_2$\\ \hline
$\om_4$ & $\om_2$ & $\om_2$ \\ \hline
\end{tabular}
\caption*{$g_4$}
\label{tab:g4}
\end{subtable}
\quad
\hspace*{-4em}
\begin{subtable}{0.3\linewidth}
\centering
\begin{tabular}{|c|c|c|c|}
\hline
$\om_1$ & $\om_2$ & $\om_1$ & $\om_2$\\ \hline
$\om_3$ & $\om_4$ & $\om_4$ & $\om_3$\\ \hline
$\om_1$ & $\om_4$ & $\om_1$ & $\om_5$\\ \hline
$\om_3$ & $\om_2$ & $\om_6$ & $\om_2$ \\ \hline
\end{tabular}
\caption*{$g_5$}
\label{tab:g5}
\end{subtable}
\quad
\hspace*{-3.5em}
\begin{subtable}{0.3\linewidth}
\centering
\begin{tabular}{|c|c|}
\hline
$\om_1$ & $\om_1$ \\ \hline
$\om_1$ & $\om_2$ \\ \hline
\end{tabular}
\caption*{$g_6$}
\label{tab:g6}
\end{subtable}

\vspace{0.5cm}

\hspace*{4em}
\begin{subtable}{0.3\linewidth}
\centering
\begin{tabular}{|c|c|}
\hline
$\om_1$ & $\om_2$ \\ \hline
$\om_2$ & $\om_1$ \\ \hline
\end{tabular}
\caption*{$g_7$}
\label{tab:g7}
\end{subtable}
\quad
\hspace*{-4em}
\begin{subtable}{0.3\linewidth}
\centering
\begin{tabular}{|c|c|c|}
\hline
$\om_1$ & $\om_1$ & $\om_2$\\ \hline
$\om_3$ & $\om_4$ & $\om_3$\\ \hline

\end{tabular}
\caption*{$g_8$}
\label{tab:g8}
\end{subtable}
\quad
\hspace*{-3.5em}
\begin{subtable}{0.3\linewidth}
\centering
\begin{tabular}{|c|c|c|}
\hline
$\om_1$ & $\om_1$ & $\om_2$\\ \hline
$\om_4$ & $\om_5$ & $\om_2$\\ \hline
$\om_4$ & $\om_3$ & $\om_3$ \\ \hline
\end{tabular}
\caption*{$g_9$}
\label{tab:g9}
\end{subtable}

\caption{Nine game forms. Forms $g_1$ - $g_6$ are tight, 
forms $g_7$ - $g_9$ are not.}
\label{9-game-forms}
\end{figure*}

Tightness is equivalent with SP-solvability 
\cite{EF70,Gur73} and with NE-solvability 
\cite{Gur75,Gur89,GN22}.  

\subsubsection{Not tight game forms}
Since tightness and SP-solvability are equivalent, 
the Shapley Theorem implies 
that all minimal not tight game forms are of size $2 \times 2$.
Two sets of outcome corresponding 
to two diagonals are disjoint.  
There are three such game forms:

\begin{figure*}[h]
%\centering
\captionsetup[subtable]{position = below}
\captionsetup[table]{position=top}
\hspace*{4em}
\begin{subtable}{0.3\linewidth}
\centering
\begin{tabular}{|c|c|}
\hline
$\om_1$ & $\om_2$ \\ \hline
$\om_2$ & $\om_1$ \\ \hline
\end{tabular}
%\caption*{$g_1$}               
\label{tab:g1_}
\end{subtable}
\hspace*{-3em}
\begin{subtable}{0.3\linewidth}
\centering
\begin{tabular}{|c|c|}
\hline
$\om_1$ & $\om_2$\\ \hline
$\om_3$ & $\om_1$\\
\hline
\end{tabular}
%\caption*{$g_2$}
\label{tab:g2_}
\end{subtable}
\hspace*{-2.5em}
\begin{subtable}{0.3\linewidth}
\centering
\begin{tabular}{|c|c|c|}
\hline
$\om_1$ & $\om_2$\\ \hline
$\om_3$ & $\om_4$\\ \hline
\end{tabular}
%\caption*{$g_3$}               
\label{tab:g3}
\end{subtable}
\caption{Three not tight $2 \times 2$ game forms}
\label{f9}
\end{figure*}

This result was strengthened in \cite{BGM09}, 
where it was shown that these three game forms 
are the only locally minimal not tight ones. 
In other words, family  $\cF$  
of not tight game forms is convex and 
class  $\cM(\cF) = \cL\cM(\cF)$  
consists of the above three game forms. 

Here we will strengthen this result further as follows. 

\begin{theorem} 
Family  $\cF$  of not tight game forms   
is strongly convex but not weakly hereditary. 
\end{theorem} 

\proof 

Obviously, a game form with a single outcome is tight. 

Consider game forms with two outcomes 
$\om_1 = a$  and  $\om_2 = b$. 
Obviously, such game form  $g$  is tight 
if and only if  one of the following cases holds: 

Case (rca): 
Game form   $g$  contains an $a$-row and an $a$-column, 
that is, there exists an $x_0 \in X$   and  $y_0 \in Y$ 
such that  $g(x,y) = a$  whenever $x = x_0$  or  $y = y_0$. 

Case (rcb): Game form  $g$  contains a $b$-row and a $b$-column. 

Case (rab): Game form $g$  contains an $a$-row and a $b$-row. 

Case (cab): Game form $g$  contains an $a$-column and a $b$-column. 

\bigskip 

\begin{example}
Family $\cF$  is not weakly hereditary,  
consider the following  $4 \times 4$  game form:

\begin{figure*}[h]
\centering
\begin{tabular}{|c|c|c|c|}
\hline
$a$ & $b$ & $a$ & $b$\\ \hline
$b$ & $a$ & $a$ & $b$\\ \hline
$a$ & $a$ & $a$ & $b$\\ \hline
$b$ & $b$ & $b$ & $a$\\ \hline
\end{tabular}
\end{figure*}

Obviously, it is not tight but becomes tight 
if we delete 
the last rows or column 
(or two last rows or columns).  
In contrast, every other, not last, 
row or column can be deleted 
and the the obtained reduced game form remains not tight. 
\end{example}

This example proves that family $\cF$ 
is not weakly hereditary, 
yet, it does not disprove strong convexity. 
Actually, family  $\cF$ is strongly convex. 
To show this, consider a not tight game form  
$g : X \times Y \rightarrow \Om$  
of size larger than  $2 \times 2$  
and fix a $2 \times 2$  not tight subform  
$g^* : X^* \times Y^* \rightarrow \Om$  in it, 
that is, $|X^*| = |Y^*| = 2$,  
$X^* \subseteq X, \; Y^* \subseteq Y$, 
and at least one of these two containments is strict.
Wlog, we can assume that  $g^*$  is formed 
by the first two rows and columns of  $g$. 

It is enough to show that one can delete 
a row $x \in X \setminus X^*$   or 
a column  $y \in Y \setminus Y^*$  such that 
the reduced game form  $g'$  is still not tight.

It is both obvious and well-known that 
merging outcomes respects tightness. 
Hence, wlog, we can assume that 
$\Om = \{\om_1,\om_2\} = \{a,b\}$ 
consists of two outcomes, and that  $g^*$  is  

\begin{figure*}[h]
\centering
\begin{tabular}{|c|c|}
\hline
$a$ & $b$ \\ \hline
$b$ & $a$\\ \hline
\end{tabular}
\end{figure*}

It is also clear that
adding an $a$-row or $a$-column 
to a tight game form respects its tightness. 
In other words, deleting an $a$-row or an $a$-column 
from a not tight game form $g$ 
respects its non-tightness.  
By symmetry, the same holds for 
$b$-rows and $b$-columns as well. 

\medskip 

Assume for contradiction that   
after deleting a row $x \in X \setminus X^*$, 
or a column from  $y \in Y \setminus Y^*$ 
from $g$, the obtained reduced subform  $g'$ is tight.  
Then, as we know, one of the cases: 
(rca), (rcb), (rab), (cab)  holds. 
Assume wlog that we delete a column 
(rather than a row) and consider all four cases. 

In case  (rca)  $g'$  
has an  $a$-column  $x \in X \setminus X^*$. 
Deleting this column from  $g$  
we obtain a not tight game form, 
which is a contradiction. 
By symmetry, case (rca) resolved too.   

Case (cab)  is trivial, 
since then both 
$g'$ and $g$  are tight, which is a contradiction.  

Thus, only case  (rab)  remains,  
and we can assume that 
after deleting each column  
$y \in Y \setminus Y^*$  from  $g$ 
the obtained reduced submatrix  $g'$  has 
both $a$- and $b$-rows. 
This is possible (only) when  $g$  contains 
4 rows and 3 columns. 

Yet, by symmetry, we can also assume that 
after deleting each row  
$x \in X \setminus X^*$  from  $g$ 
the obtained reduced submatrix  $g''$  has 
both $a$- and $b$-columns. 

This is already impossible. 
To see it, consider the $4 \times 4$  game form shown in Fig. \ref{f00} in which we can assign 
$a$ or $b$  arbitrarily to every 
symbol  $*$. 
It is not difficult to verify that 
every such assignment results in contradiction: 
the obtained subform is either tight, 
or not tight and among 
its last two rows and columns there is at least one 
deleting which results in a subform 
that is still not tight. 
\qed

\begin{figure}[h]
\centering
\begin{tabular}{|c|c|c|c|}
\hline
$a$ & $b$ & $a$ & $b$\\ \hline
$b$ & $a$ & $a$ & $b$\\ \hline
$a$ & $a$ & $*$ & $*$\\ \hline
$b$ & $b$ & $*$ & $*$\\ \hline
\end{tabular}
\caption{Case (rab), a contradiction.}
\label{f00}
\end{figure}
\noindent 

\subsubsection{Tight game forms}
Tightness is not hereditary. 
For example, game form  $g_3$  is tight 
but after deleting the last column one obtains $g_8$,  
which is not tight. 

Moreover, it was shown in \cite{BGM09}  that 
the family $\cF$ of tight game forms is not convex. 
Class $\cM(\cF)$  contains only  $1 \times 1$ game forms,  
while  $\cL\cM(\cF)$  is a complicated class, 
which seems difficult to characterize; 
only some necessary and some sufficient conditions 
are obtained in \cite{BGM09}. 
Note that  $g_3 \not\in \cL\cM(\cF)$, 
since deleting its row keeps tightness.  

\subsubsection{Totally tight game forms}
A game form is called {\em totally tight (TT)}  
if every its subform is tight;  
for example,  $g_3$  in Figure \ref{9-game-forms} is TT.  

\begin{proposition} (\cite{BGMP10}). 
Tightness of all  $2 \times 2$  subforms 
% of a game form already implies its 
already implies total tightness. 
\qed 
\end{proposition}

Sketch of the proof. 
This result is implied  
by the following two criteria of solvability  
\cite{Gur73,Sha64}. 
The first states that a game has a SP if 
(but not only if)  
every its $2 \times 2$ subgame has a SP \cite{Sha64}.  
The second claims that a game form is {\em zero-sum-solvable}  
if and only if it is tight \cite{Gur73}; 
see more details in \cite{Gur89,GN21,GN22}. 
Let us note that the second result is implicit already in \cite{EF70}. 
\qed 

\smallskip 

It is easily seen that the next three properties 
of the $2 \times 2$  game forms are equivalent: 

\medskip 
(i) tightness, \; (ii) total tightness, \; 
(iii) presence of a constant row or column. 

\medskip 

As we know, there are only three not tight 
$2 \times 2$  game forms; 
they are given in Figure \ref{f9}. % {2x2-not-tight}.

\smallskip 

Thus, family of the TT game forms 
is characterized by these three forbidden subforms. 
Hence, this family is hereditary. 

\medskip 

It is also known that 
a game form is TT if and only if 
it is {\em acyclic}, that is, for arbitrary payoffs 
of two players the obtained game has no {\em improvement cycle}; 
see \cite{BGMP10} for the proof and precise definitions. 

An explicit recursive characterization  
of the TT game forms is also given in \cite{BGMP10}.

\subsubsection{Not totally tight game forms}

As we already mentioned, a game form is 
not TT if and only if it contains at least one 
of the three not tight  $2 \times 2$  subforms 
given in Figure \ref{f9}.  
Thus, family of not TT game forms  
is weakly hereditary.  
Obviously, it is not hereditary. 
Indeed, by deleting a row or column 
from a  $2 \times 2$ game form    %%  in Figure \ref{f9} 
one obtains a $1 \times 2$  or $2 \times 1$
game form, which is tight.

\section*{Acknowledgements}
This research was supported by Russian Science Foundation, grant 20-11-20203, 

https://rscf.ru/en/project/20-11-20203/

\end{document}